\documentclass[a4paper,10pt]{article}
\usepackage[utf8]{inputenc}
\usepackage{amsmath,amssymb,amsthm}
\usepackage{graphicx,graphics,color,booktabs}
\usepackage{fancyhdr}
\usepackage{url}
\usepackage{tikz}
\newtheorem{remark}{Remark}[section]

\author{
  Bal\'{a}zs Kov\'{a}cs\footnote{e-mail: kovacs@na.uni-tuebingen.de} \\ Mathematisches Institut, Universit\"at T\"ubingen, \\ Auf der Morgenstelle 10, D-72076 T\"ubingen, Germany
}
\title{Computing arbitrary Lagrangian Eulerian maps \\ for evolving surfaces}

\newcommand\bfx{{\mathbf x}}

\newcommand\bfX{{\mathbf X}}

\newcommand\bflambda{{\boldsymbol \lambda}}

\newcommand\bfLambda{{\boldsymbol \Lambda}}
\newcommand\bfell{{\boldsymbol \ell}}

\DeclareMathOperator{\diff}{\frac{\textrm{d}}{\textrm{d}t}}

\newcommand{\Ga}{\Gamma}

\DeclareMathOperator{\laplace}{\Delta}

\newcommand{\nb}{\nabla}

\DeclareMathOperator{\pa}{\partial}
\newcommand{\R}{\mathbb{R}}

\def \t {(t)}
\def \to {\rightarrow}

\newcommand\andquad{\qquad \textnormal{ and } \qquad }

\def\AC{{\mathcal O}\hskip-2.4pt\iota}

\newcommand{\Gat}{\Gamma(t)}

\newcommand{\nbg}{\nabla_{\Gamma}}

\newcommand{\mat}{\partial^{\bullet}}
\newcommand{\Th}{\mathcal{T}_h}

\begin{document}

\maketitle

% REQUIRED
\begin{abstract}
  The good mesh quality of a discretized closed evolving surface is often compromised during time evolution. In recent years this phenomenon has been theoretically addressed in a few ways, one of them uses arbitrary Lagrangian Eulerian (ALE) maps. However, the numerical computation of such maps still remained an unsolved problem in the literature.
  An approach, using differential algebraic problems, is proposed here to numerically compute an arbitrary Lagrangian Eulerian map, which preserves the mesh properties over time. The ALE velocity is obtained from a simple spring system, based on the connectivity of the nodes in the mesh. We also consider the algorithmic question of constructing acute surface meshes.
  We present various numerical experiments illustrating the good properties of the obtained meshes and the low computational cost of the proposed approach.
  
  
  \noindent AMS: 65M50, 34A09, 35R01
\end{abstract}

\section{Introduction}
Partial differential equations on evolving surfaces with a \emph{given velocity} $v$ have been discretized using a huge variety of methods. Probably one of the most popular is the evolving surface finite element method developed by Dziuk and Elliott in \cite{DziukElliott_ESFEM}.

As it was pointed out by Dziuk and Elliott early on, in Section~7.2 of \cite{DziukElliott_ESFEM}: \textit{''A drawback of our method is the possibility of de\-ge\-ne\-rating grids. The prescribed velocity may lead to the effect, that the triangulation $\Ga_h(t)$ is distorted''}. The same issue occurs for problems in moving domains.

In recent years this problem has been theoretically addressed in a few ways:
\begin{itemize}
	\item Spatial discretizations using \emph{arbitrary Lagrangian Eulerian} (ALE) maps for moving domains have been studied many times in the literature, see for instance \cite{FormaggiaNobile1999,FormaggiaNobile2004} and \cite{Bonito2013,Bonito2013a}, and the references therein.
	\item The ALE version of the evolving surface finite element method has been proposed by Elliott and Styles in \cite{ElliottStyles_ALEnumerics}, where a better triangulation is obtained using an ALE map, i.e., allowing the nodes of the mesh to evolve with a velocity having an additional tangential component compared to the (pure normal) surface velocity.
	%This method was later analysed in \cite{ElliottVenkataraman_ALE,KovacsPower_ALE}.
	\item Elliott and Fritz \cite{ElliottFritz_DT,ElliottFritz} constructed meshes with very good properties using the quite technical \emph{DeTurck trick}.
	%  \item ... Barrett, Garcke and N\"{u}rnberg ...
	%  \item ...
\end{itemize}
%By the evolving surface ALE approach of \cite{ElliottStyles_ALEnumerics} a better triangulation is obtained by using an arbitrary Lagrangian Eulerian map, or equivalently, allowing the nodes of the mesh to evolve with a velocity which has an additional tangential component compared to the (pure normal) surface velocity.

%\medskip
%In contrast with the many papers on the theory of ALE maps for both moving domains and evolving surfaces, there are no papers, up to our knowledge, concerning \emph{how} such an ALE map should be computed.
%
%\smallskip
%We propose an algorithm here to compute an arbitrary Lagrangian Eulerian map for closed evolving surfaces.

%\medskip
We propose here to compute an arbitrary Lagrangian Eulerian map for closed evolving surfaces by integrating a differential algebraic equation (DAE) system for the nodes. We use a not necessarily tangential ALE velocity to achieve good mesh quality, while enforcing the points to stay on the surface.
% This approach is independent from any unrealistic a priori knowledge.

To our knowledge there is no such ALE algorithm available in the literature, in contrast with the many papers on the theory of numerical methods involving ALE maps for both closed evolving surfaces and moving domains.

%\medskip
\pagebreak
Many experiments have been presented in the above references, especially see \cite{ElliottStyles_ALEnumerics,ElliottVenkataraman_ALE,KovacsPower_ALE}, where smaller discretization errors have been obtained by solving evolving surface problems on ALE meshes. The ALE maps used in these experiments were unrealistic, obtained analytically from an \emph{a priori knowledge} on the surface and its evolution, using deep understanding and structure of the signed distance function (which defines the surface). No general ideas on the computation of ALE maps for evolving surfaces have been proposed in these papers. Numerical analysis of the ALE evolving surface finite element method has been studied in \cite{ElliottVenkataraman_ALE} and \cite{KovacsPower_ALE}.

For non-ALE methods it is usual that nodes disappear and/or new nodes are added to the mesh (see, e.g., \cite[Figure~14]{ElliottFritz_DT}). From a theoretical point of view this seems a minor issue, however from the implementation side this is undesirable, since then a map -- between the old and the new mesh -- has to be computed after each remeshing process. In most cases this is a nontrivial and costly task, and hence, should be avoided.

\medskip
In the present paper we propose a general algorithm to compute a suitable ALE map, \emph{without any a priori knowledge}, for meshes of closed evolving surfaces.

The approach is based on the following idea:
Usually, the surface evolution is given by an ODE system with a surface velocity assumed to be normal. We use an additional tangential velocity for a possibly degenerated mesh to improve grid quality. In general such tangential velocities are not straightforward to construct, therefore we use a not necessarily tangential ALE velocity and introduce a constraint to keep the nodes on the surface. Altogether this is finally formulated as a DAE system for the nodes.
The ALE velocity is based on a spring system, where the nodes are connected along the edges by springs. The numerical solution of the DAE system gives the new mesh. We use implicit Runge--Kutta methods (in particular Radau IIA methods), and a more efficient splitting method, combined with explicit Runge--Kutta methods, to integrate the system in time.
%Naturally, the idea of using differential algebraic problems is independent from the particular ALE velocity proposed here.

The computation of the arbitrary Lagrangian Eulerian mesh here is free of any \emph{a priori knowledge} in the following sense: the algorithm uses the distance function at each time, but it does not use its structure or any other special properties of it, unlike the ALE maps from the literature. 
%We note here that the use of the distance function can be also avoided.

This approach for closed evolving surfaces can be used as a tool in the computation of ALE meshes for \emph{moving domains}: In \cite[Section~2.4]{FormaggiaNobile1999} arbitrary Lagrangian Eulerian maps for moving domains are obtained by solving a parabolic problem, or the corresponding stationary problem, while in \cite{FarhatLesoinneMaman} an elastodynamic equation system is used for the same purpose. However, for these approaches the evolution of the boundary still needs to be \emph{known a priori}. The problem of numerically finding such a boundary evolution has not been solved in these papers. In fact, to determine such a boundary evolution is equivalent to finding an ALE map for a lower dimensional closed surface, which is the same problem as we consider in this paper. Hence, the algorithm proposed here can also serve as a tool to compute boundary evolutions, which can be used together with the well understood classical ALE methods for moving domains, for instance the ones proposed in \cite{FarhatLesoinneMaman,FormaggiaNobile1999}.

We give some further details on possible extensions of the proposed algorithm: to handle other mesh properties (e.g., acuteness), an adaptive version, and a local version as well.

We present various numerical experiments illustrating the validity of the differential algebraic model, and also the performance of the proposed algorithm compared to the ALE maps given in the literature. We also report on computational times in the case of evolving surface PDEs.

\section{Evolving surfaces, ALE maps and PDEs}
\label{section: ALE maps}

As our main motivation lies in the numerical solution of parabolic PDEs on evolving surfaces we shortly recap the setting of \cite{DziukElliott_ESFEM}. We will also use this setting as an illustrative background to the proposed algorithm.

Let $\Ga\t \subset \R^{m+1}$, $0 \leq t \leq T$, be a smooth evolving closed hypersurface. %, with $m\leq 3$
Further, let the evolution of the surface given by the smooth velocity $v$, assumed to be normal. Let $\mat u = \pa_{t} u + v \cdot \nb u$ denote the material derivative of $u$, the tangential gradient is denoted by $\nbg$ and given by $\nbg u = \nb u - \nb u \cdot \nu \nu$, with unit outward normal $\nu$. We denote by $\laplace_\Ga = \nbg\cdot\nbg$ the Laplace--Beltrami operator.

We consider the following linear evolving surface PDE:
\begin{equation}
\label{eq: evolving surface PDE}
\begin{aligned}
\mat u + u \nb_{\Gat} \cdot v - \laplace_{\Gat} u =&\ f  & \textrm{ on } & \Ga\t ,\\
u(\cdot,0) =&\  u_0 & \textrm{ on } & \Ga(0).
\end{aligned}
\end{equation}
Basic and detailed references on evolving surface PDEs and on the evolving surface finite element method (ESFEM) are \cite{DziukElliott_ESFEM,DziukElliott_acta,DziukElliott_L2} and \cite{diss_Mansour}.
The ALE version of the evolving surface finite element method has been proposed in the paper \cite{ElliottStyles_ALEnumerics}, which also contains a detailed description and many experiments. While numerical analysis of full discretizations can be found in \cite{ElliottVenkataraman_ALE} and \cite{KovacsPower_ALE}, the former is more concerned about spatial discretizations and BDF methods of order 1 and 2, while the latter is more focused on high-order BDF and Runge--Kutta time discretizations.

\bigskip
As we aim at the numerical solution of the PDE \eqref{eq: evolving surface PDE} using the surface finite elements developed by Dziuk and Elliott \cite{DziukElliott_ESFEM}, the surface is discretized using a triangular mesh. The description of representation of evolving surfaces and that of discrete surfaces can be found in the following subsections.

\subsection{Surface representations}
\label{section: surface representation}

Evolving surfaces are usually described in two ways, which have different advantages, hence we will use both of them for various purposes.

{\bf Distance function representation.}
Based on a signed distance function the evolving $m$-dimensional closed surface $\Gamma(t)\subset\R^{m+1}$ is given by
\begin{equation*}
\Gat = \{ x\in \R^{m+1} \mid d(x,t)=0 \} ,
\end{equation*}
with a function $d:\R^{m+1} \times [0,T] \to \R$ (whose regularity depends on the smoothness of the surface), cf.~\cite{Dziuk88,DziukElliott_ESFEM}.

{\bf Diffeomorphic parametrization.} The surface can also be described by a diffeomorphic parametrization, cf.\ \cite{soldriven} and \cite{DziukElliott_ESFEM}.

We consider the evolving $m$-dimensional closed surface $\Gamma(t)\subset\R^{m+1}$ as the image
$$
\Ga(t) = \{ X(p,t) \mid p \in \Ga(0) \}
$$
of a sufficiently regular vector-valued function $X:\Ga(0)\times [0,T]\to \R^{m+1}$, where $\Ga(0)=\Ga^0$ is the smooth closed initial surface, and $X(p,0)=p$.
It is convenient to think of $X(p,t)$ as the position at time $t$ of a moving particle with label $p$, and of $\Ga(t)$ as a collection of such particles.
The parametrisation also satisfies the ODE system, for a point $p\in\Ga(0)$,
\begin{equation}
\label{eq: surface velocity}
\partial_t X(p,t)= v(X(p,t),t),
\end{equation}
where $v(\cdot,t)\in\R^{m+1}$ is the velocity of the surface.
Note that for a known velocity field $v$, the position $X(p,t)$ at time $t$ of the particle with label $p$ is usually obtained by solving the ordinary differential equation \eqref{eq: surface velocity} from $0$ to $t$ for a fixed $p$.

\subsection{Surface approximation}
\label{section: surface approximation}

The smooth initial surface $\Ga(0)$ is approximated by a triangulated surface $\Ga_h(0)$, i.e., an admissible family of triangulations $\Th(0)$ of maximal element diameter $h$; see \cite{DziukElliott_ESFEM} for the notion of an admissible triangulation, which includes quasi-uniformity.
Let $x_j(0)$, ($j=1,2,\ldots, N$) denote the nodes of $\Ga_h(0)$ lying on the initial smooth surface $\Ga(0)$. The nodes will be evolved in time with the given normal velocity $v$, by solving the ODE
\begin{equation}
\label{eq: ODE for nodes}
\diff x_j\t =v(x_j\t,t) \qquad (j=1,2,\dotsc,N),
\end{equation}
which is simply \eqref{eq: surface velocity} for the nodes $(x_j)$.
Obviously, the nodes remain on the surface $\Gat$ for all times, i.e., $d(x_j\t,t)=0$ for $j=1,2,\dotsc,N$ and for all $t\in[0,T]$.

Therefore, the smooth surface $\Gat$ is also approximated by a discrete surface $\Ga_h\t$, whose elements also form a triangulation $\Th\t$. We have
\begin{equation*}
\Ga_h(t) = \bigcup_{E(t)\in \Th(t)} E(t) .
\end{equation*}

As already mentioned in the introduction, we assume that the surface does not develop topological changes due to the evolution. The assumption on quasi-uniformity over time, i.e., there is a fixed $c>0$ (independent of $t$) such that for any triangle $E(t)\in\Th(t)$ the radius of the inscribed circle $\sigma_{E(t)}$ satisfies
\begin{equation*}
\frac{h_{E(t)}}{\sigma_{E(t)}} \leq c , %  \qquad ( t\in[0,T] )
\end{equation*}
for all $t\in[0,T]$, is generally not always satisfied during time evolution.

\begin{figure}[ht!]
	\centering
	\includegraphics[width=.32\textwidth]{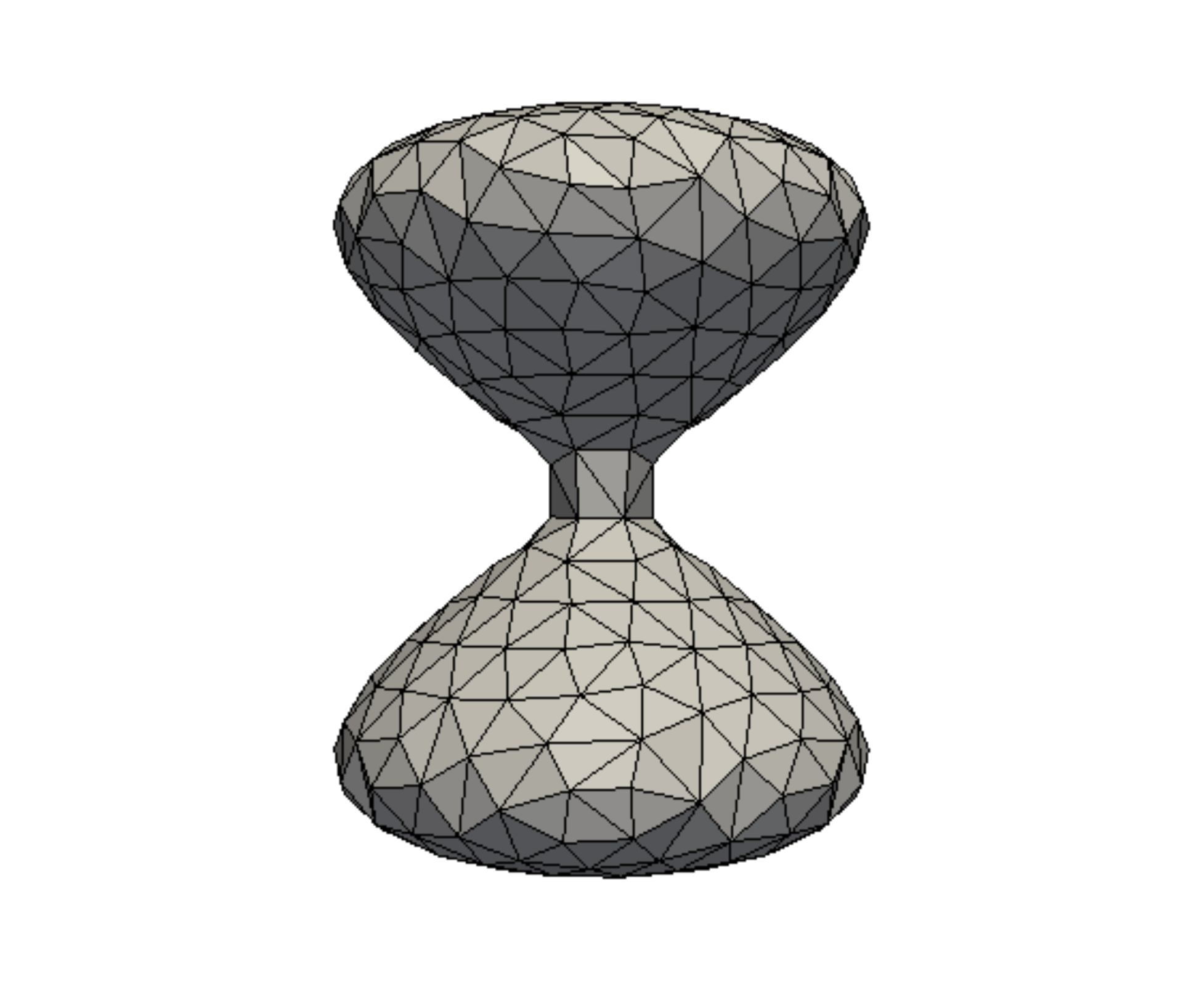}
	\includegraphics[width=.32\textwidth]{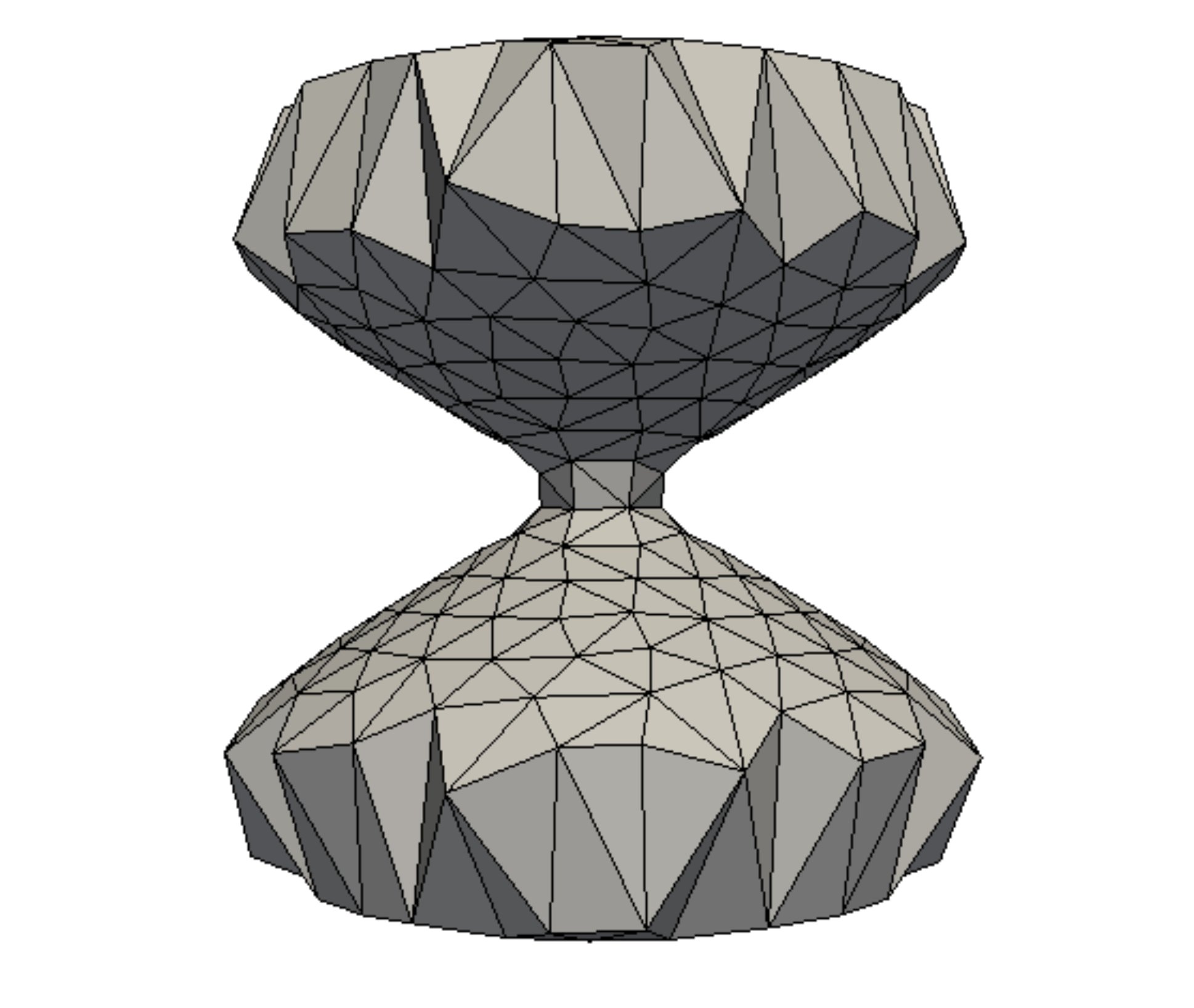}
	\includegraphics[width=.32\textwidth]{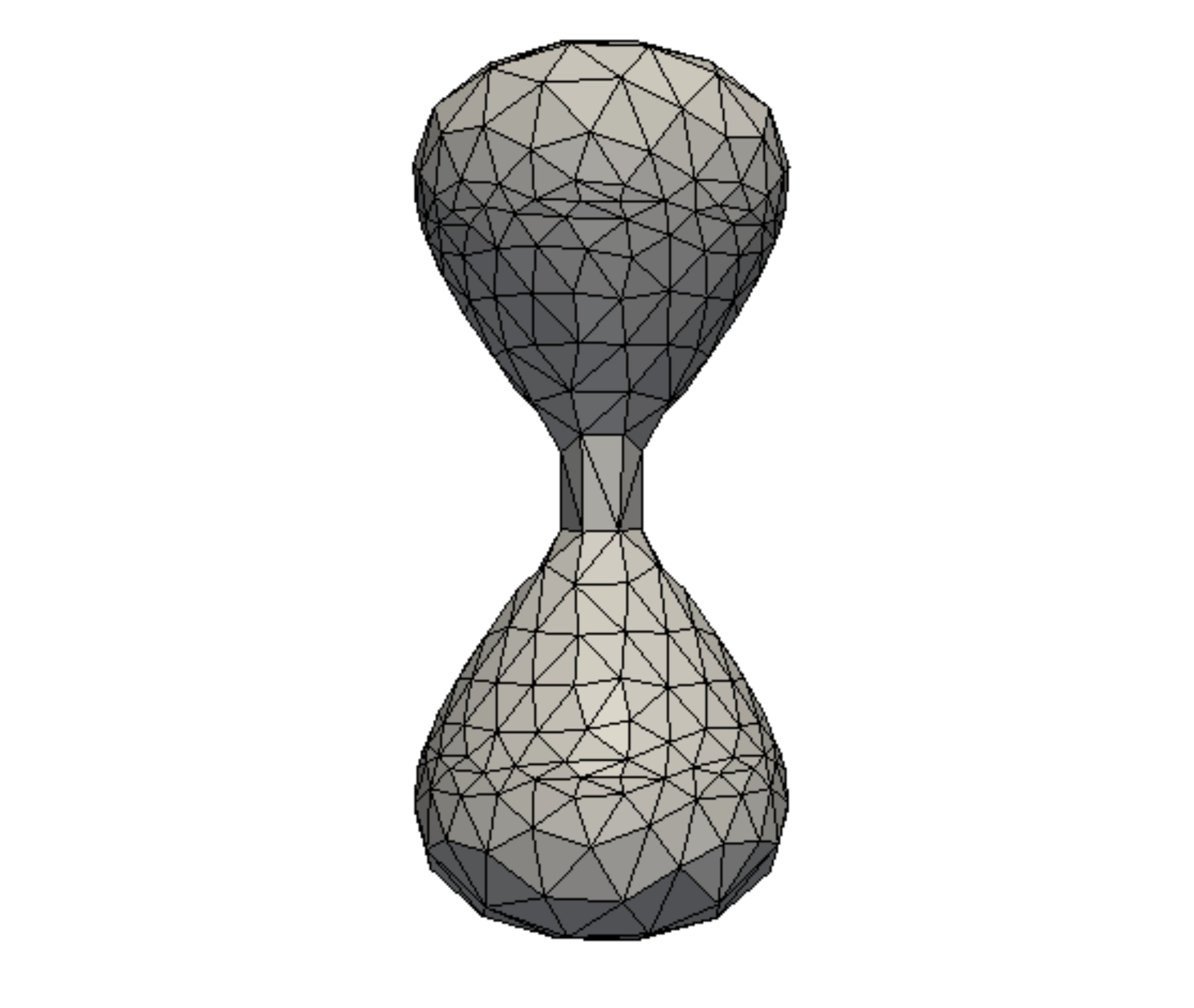}\\
	\caption{Normal evolution of a closed surface at time $t=0, 0.2, 0.6$; see also in \cite{ElliottStyles_ALEnumerics}}
	\label{figure: example}
\end{figure}

As an example, from \cite{ElliottStyles_ALEnumerics}, to de\-ge\-ne\-rating surface evolution, we evolved a surface using the ODE \eqref{eq: ODE for nodes}. As observed in Figure~\ref{figure: example}: however the initial mesh (left) is quasi-uniform, the meshes at later times (middle and right) do not preserve the good mesh qualities. Both quite bad surface resolution and unnecessarily fine elements occur.

\section{Computing ALE maps as general constrained problems}
We now propose an approach which will be used to determine a suitable ALE map for evolving surfaces, to maintain mesh quality during the surface evolution. In fact, we directly compute the new positions of the nodes.

We consider an other parametrization of $\Gamma(t)$, with good mesh properties, and which is different from $X$ in \eqref{eq: surface velocity}, called an \emph{arbitrary Lagrangian Eulerian map}. The corresponding ODE system is
\begin{equation}
\label{eq: tangential ALE velocity - analytic}
\partial_t X(p,t)= v(X(p,t),t) + w(X(p,t),t),
\end{equation}
with the pure tangential velocity $w$. 

For evolving surface problems the surface velocity $v$ is usually assumed to be known and normal. However, such pure tangential ALE velocities are not given, and also not easy to obtain, in general.
Therefore, we allow velocities (still denoted by) $w$ which improve mesh quality, but have small non-tangential components, hence may lead points away from the surface. To compensate this, a constraint is introduced in order to keep the surface $\Ga\t$ unaltered.
Therefore, the set of ODEs \eqref{eq: tangential ALE velocity - analytic} is modified into the following differential algebraic equation (DAE) system (of index 2) with Lagrange multiplier $\lambda$, for $p\in\Ga(0)$,
\begin{equation}
\label{eq: constrained ALE - analytic}
\begin{aligned}
\partial_t X(p,t) =&\ v(X(p,t),t) + w(X;X(p,t),t) - D(X(p,t))^T \lambda(p,t) , \\
d(X(p,t),t) =&\ 0 ,
\end{aligned}
\end{equation}
where $D(X) = \pa d(\cdot,t)/\partial X$. The first argument $X$ in $w$ indicates that the additional velocity may (and usually does) depend on the whole surface. As it causes no confusion we will drop this argument later on.

The analogous differential algebraic equation system for the nodes of the surface approximation mesh $\Ga_h\t$, which are collected into the vector $\bfx\t=(x_j\t)_{j=1}^N$, and with Lagrange multiplier $\bflambda(t) \in \R^N$, reads as
\begin{equation}
\label{eq: constrained ALE ODE system}
\begin{aligned}
\diff \bfx\t =&\  v(\bfx\t,t) + w(\bfx\t,t) - D(\bfx\t)^T \bflambda\t , \\
d(\bfx\t,t) =&\ 0,
\end{aligned}
\end{equation}
with initial value $\bfx(0)=(x_j(0))_{j=1}^N$. The constraint $d(\bfx\t,t) = 0$ is meant pointwise, i.e., as $d(x_j\t,t) = 0$ for $j=1,2,\dotsc,N$, while the matrix $D(\bfx) = \pa d(\cdot,t) / \pa \bfx$. 
Concerning notation: we will apply the convention to use boldface letters to denote vectors in $\R^{3N}$ or $\R^N$ collecting nodal values of discretized variables.

Naturally, the DAE system is independent of the choice of the ALE velocity $w$, if there is some (physical, biological or modeling) knowledge on the general type of the surface evolution a user can propose a suitable ALE velocity accordingly.

In the next sections we will propose a very intuitive way to define the ALE velocity $w$, and we will also discuss numerical methods for the solution of the DAE system.

%\subsection{Laplacian smoothing}
%
%...

\subsection{A spring system based arbitrary ALE velocity}
We use here a simple idea to determine the velocity $w$: let us assume that the nodes of the mesh are connected by springs following the edges of the elements, i.e., the topology of the spring system is determined by the triangulation $\Ga_h\t$.

% Our main idea is that, starting from a high-quality triangulation, due to pure normal surface evolution (with usually small time stepsizes) the mesh is not getting too far from equilibria (e.g., as observed in Figure~\ref{figure: example}), and that in the equilibrium state most of the springs (edges) have almost the same length, and the triangles are close to regular triangles. We propose definitions to $F$, in such a way that initially the springs over long edges are stretched, short ones are compressed, while springs with average lengths are in, or close to, equilibrium.

\begin{figure}[htbp]
	\centering
	\begin{tikzpicture}
	\filldraw
	(1,1) circle (2pt) node[align=left,   above left] {$(x_j)^e$} --
	%(1.65,1) node[align=left,   below] {$\nu_e$} --
	(2.6,1) node[align=left,   above] {$e$} --
	(4,1) circle (2pt) node[align=right,  above] {\; $x_j$};
	\draw
	(4,1) -- (3.7,1.9) (4,1) -- (5.1,1.71) node[align=right, above right]{$\{(x_j , \cdot )\}$} (4,1) -- (4.9,0.3) (4,1) -- (3.65,0.15);
	\draw[->,very thick]
	(4,1) -- (5,1) node[align=left,   below left] {$\nu_e$};
	\end{tikzpicture}
	\caption{A typical node $x_j$ on an edge $e$, the set of edges $\{(x_j , \cdot )\}$, etc.\ used to define $F_j(\bfx)$}
	\label{fig: edge defs}
\end{figure}
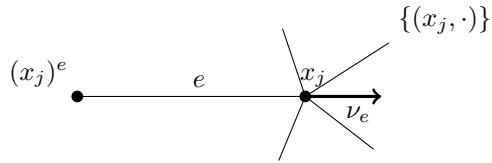

This system defines a force function $F$, which we use to define the ALE velocity by setting 
\begin{equation}
\label{eq: spring based w}
w(\bfx,t) = k F(\bfx), \qquad \textnormal{with a spring constant $k$ chosen later}. 
\end{equation}
The force function $F$ is computed based on the connectivity (described by the elements), and by the forces over the edges based on a length function $\ell_p$ (the desired length of springs). The net force $F_j(\bfx)$ at a node $x_j$ is given by, for $j=1,2,\dotsc,N$,
\begin{equation}
\label{eq: force function}
F_j(\bfx) = \sum_{e\in \{(x_j , \cdot )\}} f(e),
\end{equation}
where the set $\{(x_j , \cdot )\}$ collects all the edges $e=(x_j , (x_j)^e)$ having $x_j$ as one of their nodes, while $(x_j)^e$ simply denotes the other node across the edge $e$, see Figure~\ref{fig: edge defs}. Then $f(e)$ is the force along the edge $e$, given by
\begin{equation*}
f(e) =  \big(\ell_p(e) - |e|\big) \ \nu_e, \qquad \textnormal{ with unit vector } \quad \nu_e = \frac{x_j - (x_j)^e}{\|x_j - (x_j)^e\|} , %\quad \textnormal{ and } \k=1 ,
\end{equation*}
and current length $|e|$.

We propose to use the following length function:
\begin{equation}
\label{eq: length function}
\begin{aligned}
\ell_p(e) =&\
\begin{cases}
m_e + (1-p)(M_e - m_e), & \textnormal{ if } |e| \geq m_e + (1-p)(M_e - m_e) , \\
m_e + p(M_e - m_e), & \textnormal{ if } |e| \leq m_e + p(M_e - m_e) , \\
|e|, & \textnormal{ otherwise} ,
\end{cases} \\
&\ \qquad \qquad \qquad \qquad \hbox{where } \quad (m_e=\min_e |e|, \ M_e=\max_e |e|, \ p\in(0,1)).
\end{aligned}
\end{equation}
%where $m_e=\min_e |e|$,  $M_e=\max_e |e|$, $p\in(0,1))$.

Other length functions can also be used, for instance the $C^k$ version of $\ell_p$, or more complicated functions, such as the analogue of $\ell_p$ with more steps. Force functions based on inverse edge length are also possible to use. 

The same force function \eqref{eq: force function} has been used by Strang and Persson in DistMesh \cite{distmesh}. However, constrained systems are not appearing there, and our length function has significant differences compared to the one used in DistMesh. They have not used their approach to compute ALE maps. In particular compared to DistMesh, we do not add or delete nodes (which is essential in DistMesh for the meshing), furthermore we compute the ideal spring length in a different way (DistMesh having a rule which fits better with the possibility to delete nodes and also with the startup of their process). 

As Strang and Persson encourage their users to \textit{''[...] modify the code according to their needs''}, we indeed adapt DistMesh to suit our purposes: the force function is computed using their modified code.

\subsection{The DAE system}

By plugging in the velocity rule \eqref{eq: spring based w} into the DAE system \eqref{eq: constrained ALE ODE system} we obtain
\begin{equation}
\label{eq: constrained spring system}
\begin{aligned}
\diff \bfx\t =&\  v(\bfx\t,t) + k F(\bfx\t) - D(\bfx\t)^T \bflambda\t , \\
d(\bfx\t,t) =&\ 0 .
\end{aligned}
\end{equation}
with the spring constant $k$, which has to be chosen by the user based on the problem at hand. The numerical approximation of this system immediately gives the position of the nodes under the ALE map.

In the following we will propose some numerical methods for the time integration.

\subsubsection{Runge--Kutta solution of the DAE system}

For a given stepsize $\tau>0$, an $s$-stage implicit Runge--Kutta method, see \cite[Section~VII.4]{HairerWannerII}, applied to the DAE system \eqref{eq: constrained spring system} determines solution approximations $\bfx^{n+1}$ and approximations to the Lagrange multiplier $\bflambda^{n+1}$, as well as the internal stages $\bfX_{ni}, \bfLambda_{ni}$ by the equations
\begin{align}
\label{eq: RK internal stage a}
\bfX_{ni} &= \bfx^n + \tau \sum_{j=1}^s a_{ij} {\cal F}(\bfX_{nj},\bfLambda_{nj}), && i=1,2,\ldots,s, \\
\label{eq: RK internal stage lambda}
\bfLambda_{ni} &= \bflambda^n + \tau \sum_{j=1}^s a_{ij} \bfell_{nj}, && i=1,2,\ldots,s, \\
\label{eq: RK sol}
\bfx^{n+1} &= \bfx^n + \tau \sum_{i=1}^s b_i {\cal F}(\bfX_{nj},\bfLambda_{nj}), 
& \textnormal{and } \quad &
\bflambda^{n+1} = \bflambda^n + \tau \sum_{i=1}^s b_i \bfell_{nj},
\end{align}
where, for $i=1,2,\ldots,s$,
\begin{align}
\label{eq: RK eq}
{\cal F}(\bfX_{ni},\bfLambda_{ni}) =&\ v(\bfX_{ni},t_{ni}) + k F(\bfX_{ni}) - D(\bfX_{ni})^T \bfLambda_{ni} ,  \\
\label{eq: RK constraint}
\textnormal{and } \qquad d(\bfX_{ni},t_{ni})=&\ 0 .
\end{align}
First the system \eqref{eq: RK internal stage a}, \eqref{eq: RK eq}, \eqref{eq: RK constraint} is solved, using a simplified Newton iteration, then one can compute $\bfell_{nj}$ from \eqref{eq: RK internal stage lambda}, and finally obtain $\bfx^{n+1}$, $\bflambda^{n+1}$ from \eqref{eq: RK sol}.

The method is determined by its Butcher-tableau, i.e., the coefficient matrix $\AC = (a_{ij})_{i,j=1}^s$ and its vector of weights $b=(b_i)_{i=1}^s$, and $c_i=\sum_{j=1}^s a_{ij}$ ($i=1,2,\dotsc,s$).

In the following we will use the Radau IIA methods, which are $s$-stage Runge--Kutta methods of classical order $2s-1$.
More details on index 2 DAEs in general, as well as on their Runge--Kutta approximations can be found in \cite{HairerLubichRoche}, or \cite[Chapter~VII.]{HairerWannerII}.

\subsubsection{Splitting methods for the constrained system}
\label{section: splitting}
We propose to solve the DAE \eqref{eq: constrained spring system} using a splitting method, in order to decrease computational time. This is supported by making the following observations: 
The implicit Runge--Kutta solution of the DAE system is usually expensive.
The normal movement is a pointwise nonstiff ODE system, without constraint (if integrated exactly), hence usually can be solved very cheaply.
Also there are some examples where the normal velocity is not explicitly given, the surface evolution is defined by a direct mapping of a reference surface, i.e., we cannot use the system \eqref{eq: constrained spring system}, see for example the surface \cite[Test~Problem~2]{DziukKronerMuller}, \cite[equation~(6.1)]{LubichMansourVenkataraman_bdsurf}.

The most straightforward way is to split the differential algebraic problem \eqref{eq: constrained spring system} in a way that one of the subproblems is the original ODE system \eqref{eq: ODE for nodes}. The main benefit is that this subproblem automatically satisfies the constraint up to a small error.

The Lie splitted DAE system, over the time interval $[t_n,t_{n+1}]$, reads as
\begin{equation}
\label{eq: splitted system}
\left\{
\begin{aligned}
\diff \bfx^v(t) =&\ v(\bfx^v(t),t) , \\
&\ \\
\bfx^v(t_n) =&\ \bfx(t_n) ,
\end{aligned}
\right.
\qquad \quad
\left\{
\begin{aligned}
\diff \bfx^w(t) =&\ k F(\bfx^w(t)) - D(\bfx^w(t))^T \bflambda\t , \\
d(\bfx^w(t),t_{n+1}) =&\ 0, \\
\bfx^w(t_n) =&\ \bfx^v(t_{n+1}),
\end{aligned}
\right.
\end{equation}
and finally, setting $\bfx(t_{n+1}) = \bfx^w(t_{n+1})$, and where $F$ is a function is defined in \eqref{eq: force function}, while the constraint is still meant coordinatewise and $D(\bfx) = \pa d(\cdot,t_{n+1}) / \pa \bfx$. In fact, this order of subproblems allows us to drop the constraint from the $v$-system, as the second system will eliminate the small errors mentioned above.

In order to decrease computational time, we solve the $v$-system using the explicit Euler method. We do not require a high-order method here, as the mesh quality depends only on the second step.
For the time integration of the $w$-system we use the classical $4$-stage Runge--Kutta method and then projecting back to the surface. Due to the stiffness of the problem, it is approximated through several substeps over interval $[t_n,t_{n+1}]$.
We choose this method, since being explicit, it is very cheap, having only four right-hand side evaluations per step, but still have high-order.

We make some general remarks on the methods described above.
\begin{remark}
	If a parabolic PDE is numerically solved on the discrete surface $\Ga_h(t_{n})$ then the matrices are assembled on the improved mesh, using the ALE scheme, cf.\ \cite[(4.3) and (4.4)]{ElliottVenkataraman_ALE} or \cite[(3.5) and (3.6)]{KovacsPower_ALE}.
	The tangential vector appearing in the formulas for ALE ESFEM can be computed from the obtained nodes and from the normal velocity.
	
	The algorithm clearly uses the surface velocity $v$ in the exact same way as for the purely Lagrangian ESFEM algorithms: only to evolve the surface. However, the distance function is used as the constraint (being an inexpensive pointwise operation), which is avoided by the standard ESFEM method.
	
	It is usual to use the same time discretization method for the surface evolution as for the time integration of the parabolic PDE on the surface. This time discretization scheme could be also used to solve the DAE or the splitted system.
\end{remark}

\begin{remark}
	Higher order splitting methods, or splittings where the order of subproblems is reversed would not improve the mesh quality. Since, the mesh quality only depends on the sufficiently good numerical solution of $w$-system, while the $v$-system takes care of the surface evolution. For example, a reversed Lie splitting would first construct a good mesh, then evolve it over time in the normal direction, which can still lead to mesh distortions.
\end{remark}

\section{Possible extensions}

We now briefly describe some further extensions to the above approach. They are rather straightforward to implement, except the one using an approximative distance function (such approximative distance functions are studied, e.g., in \cite{Chang,Hoppe} and the references therein). However, thoroughly comparing and reporting on these extensions would expand the paper enormously, therefore we will here restrict our numerical experiments to the case described above.

\subsection{Construction of acute or nonobtuse surface meshes}
\label{section: acute meshes}

In a couple of recent works discrete maximum principles (DMPs) and invariant regions have been studied for surface PDEs discretized by surface finite elements, see \cite{lumped} and \cite{surfDMP}. It is well known that \emph{acute} or \emph{nonobtuse} meshes are required for DMPs even for flat domains, and also for triangulated surface meshes. In general the meshes generated by usual algorithms (for instance DistMesh \cite{distmesh}, the grid generators of DUNE \cite{DUNE}, etc.) do not necessarily satisfy these angle conditions.

The force function proposed here can be modified in such a way that the resulting mesh is acute or nonobtuse during evolution. In practice, the function $F$ should be obtained not only based on the length function \eqref{eq: length function}, but also on an \emph{angle function}, which does not allow angles approaching a prescribed tolerance $\alpha_{\textnormal{TOL}}$ given by the user. This can be viewed as three cords added to a triangle at each angle, which does not allow the angles to expand above the specified angle $\alpha_{\textnormal{TOL}}$. 

\subsection{Approximating the distance function}

The usage of the distance function could also be completely avoided in the splitting scheme. By using a modified distance function $\widetilde d(\cdot,t_{n+1})$, which is an approximation of $d(\cdot,t_{n+1})$ based on the nodes obtained from the $v$-system. Such an approximation can be obtained by a suitable interpolation process over the nodes $\bfx^v(t_{n+1})$.
For instance by a spline interpolation over the nodes, this problem appears to be well studied in the computer science literature, see for example \cite{Hoppe,Chang} and the references therein.

This approach can be useful also in cases where the surface evolution is coupled to the problem on the surface, such as \cite{DziukMCF}, or \cite{soldriven}.

We note here that in the case of such a modified algorithm the nodes $\bfx(t_{n+1})$ no longer stay on the exact surface for all times. Therefore, the spatial convergence results of \cite{ElliottVenkataraman_ALE,KovacsPower_ALE} do not hold. However, in order to show convergence results the approach of \cite{soldriven} -- using three surfaces: the exact surface, its interpolation and the discrete surface -- can be used in this case.

\subsection{Adaptivity}

A cheap adaptive method can be obtained by using some mesh quality test during the time integration of the \emph{autonomous} $w$-system of the splitting approach. In this case the following algorithm could be realized -- assuming a fast implementation of the mesh quality tester. First test the mesh quality, if it is good accept it; else perform a timestep solving the $w$-system in \eqref{eq: splitted system}, and then repeat. In such a case the numerical solution of the DAE system is avoided for already good meshes, whereas the computational overhead due to the mesh quality test is negligable.

Later on we will give some possible mesh quality measures, but it can be anything suitable specified by the user, cf.~\cite{Field_meshqualy}.

\subsection{A local version of the constrained problem}

To further decrease computational cost of the Runge--Kutta or the splitting method one may use an ALE map only locally. Since in most cases the evolution of the discrete surface distorts the mesh only locally, one could integrate the constrained system only on those patches which consist of ill shaped triangles, (e.g., those with too small or large angles), and additionally a few layers of neighbouring elements. The rest of the nodes are only evolved by the surface velocity.

\section{Numerical experiments}

\subsection{ALE map tests}
%Possible examples:
%\begin{itemize}
%  \item Elliott \& Styles example (as benchmark example);
%  \item Elliott \& Venkataraman example: surface with four holes (Example 6.3) (possible to use as benchmark example);
%  \item Elliott \& Venkataraman example: surface with changing conormal example (Example 6.4) (possible to use as benchmark example);
%\end{itemize}

%\noindent Possible measures and plots:
%\begin{itemize}
%  \item compare with ''perfect ALE maps'';
%  \item CPU times of DistMesh steps per time level;
%  \item look for mesh quality indicators;
%  \item run some parabolic PDEs;
%\end{itemize}

Now we will present some numerical experiments validating the choice of the ALE velocity \eqref{eq: spring based w} based on the spring system, and also illustrating the good qualities of the DAE model \eqref{eq: constrained spring system}. We also compared the Runge--Kutta and the splitting approach to the pure normal evolution of the surface and also to the ALE maps given in the literature. These examples have been used many times previously, see for instance \cite{ElliottStyles_ALEnumerics,ElliottVenkataraman_ALE,KovacsPower_ALE} and the references therein.
Through these ALE maps, which we call literature ALE maps, it will be also clarified further what was meant under \emph{a priori knowledge} in the introduction.

%During the time integration of \eqref{eq: splitted system} we use the classical 4-stage Runge--Kutta method with $25$ substeps, and $p=0.4$ in \eqref{eq: length function}.

\subsubsection{A dumbbell-shaped surface \cite{ElliottStyles_ALEnumerics}}
\label{section: E-S example}
Let the closed surface $\Ga\t$ be given by the zero level set of the distance function
\begin{equation}
\label{eq: surface E-S}
d(x,t) = x_1^2 + x_2^2 + K\t^2 G\Big( \frac{x_3^2}{L\t^2} \Big)-K\t^2, \quad \textrm{i.e.,} \quad
\Ga\t = \{x\in \R^3 \ \big| \ d(x,t)=0\}.
\end{equation}
Here the functions $G$, $L$ and $K$ are given by
\begin{align*}
G(s) =&\ 200 s \Big( s - \frac{199}{200}\Big) , \\
L(t) =&\ 1 + 0.2\sin(4\pi  t) , \\
K(t) =&\ 0.1 + 0.05\sin(2\pi  t) .
\end{align*}
The normal velocity $v$ describes the surface evolution at the nodes by the ODEs:
\begin{equation}
\label{eq: normal movement ODE E-S}
\diff x_j\t = v(x_j\t,t)
\end{equation}
for $j=1,2,\dotsc,N$. The surface velocity $v$ in $x_j\t$ is given by
\begin{equation}
\label{eq: v def by distance function}
v(x_j\t,t) = V_j \nu_j, \quad \textnormal{where} \quad V_j=\frac{-\pa_t d(x_j\t,t)}{|\nb d(x_j\t,t)|}, \quad \nu_j=\frac{\nb d(x_j\t,t)}{|\nb d(x_j\t,t)| } .
\end{equation}

Finally, the literature ALE map is given by
\begin{equation}
\label{eq: ALE map E-S}
(x_i\t)_1= (x_i(0))_1 \frac{K\t}{K(0)}, \quad (x_i\t)_2= (x_i(0))_2 \frac{K\t}{K(0)}, \quad (x_i\t)_3= (x_i(0))_3 \frac{L\t}{L(0)},
\end{equation}
for every $t \in [0,T]$ and for $i=1,2,\dotsc,N$, as suggested in \cite{ElliottStyles_ALEnumerics}. This map clearly uses a priori knowledge on the structure of the distance function \eqref{eq: surface E-S}.

\begin{figure}[ht!]
	\centering
	\includegraphics[width=.99\textwidth]{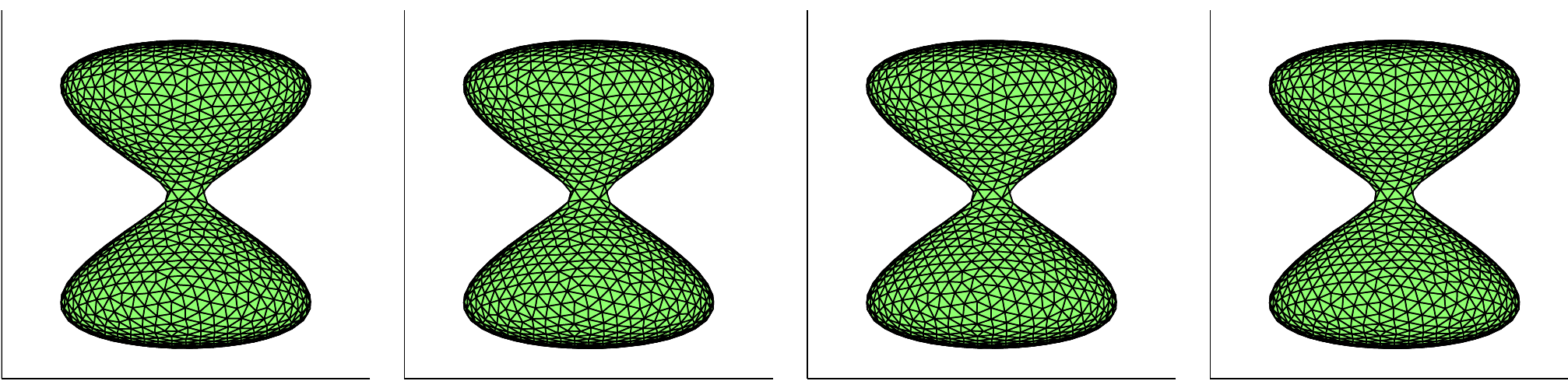}\\
	\includegraphics[width=.99\textwidth]{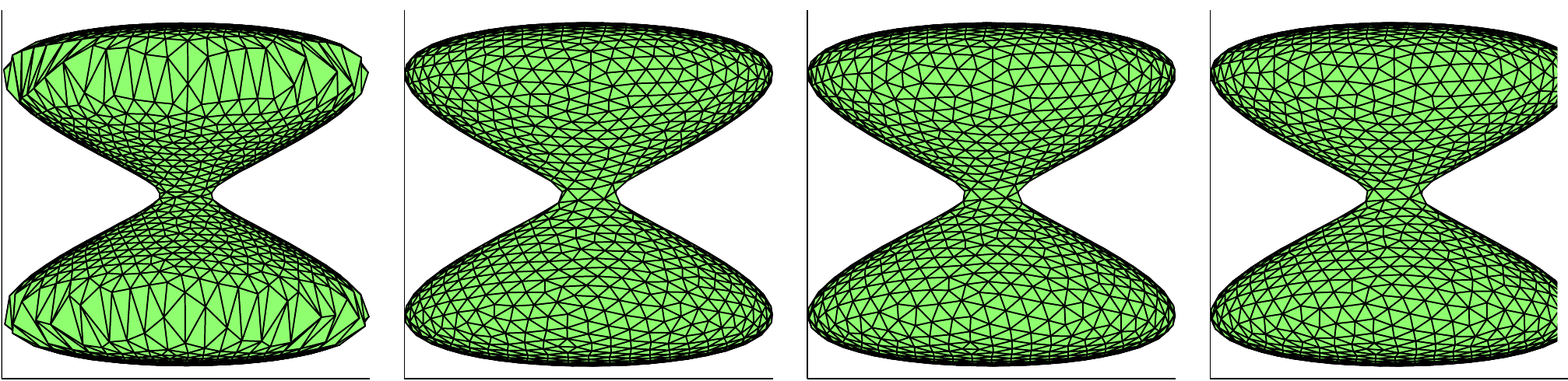}\\
	\includegraphics[width=.99\textwidth]{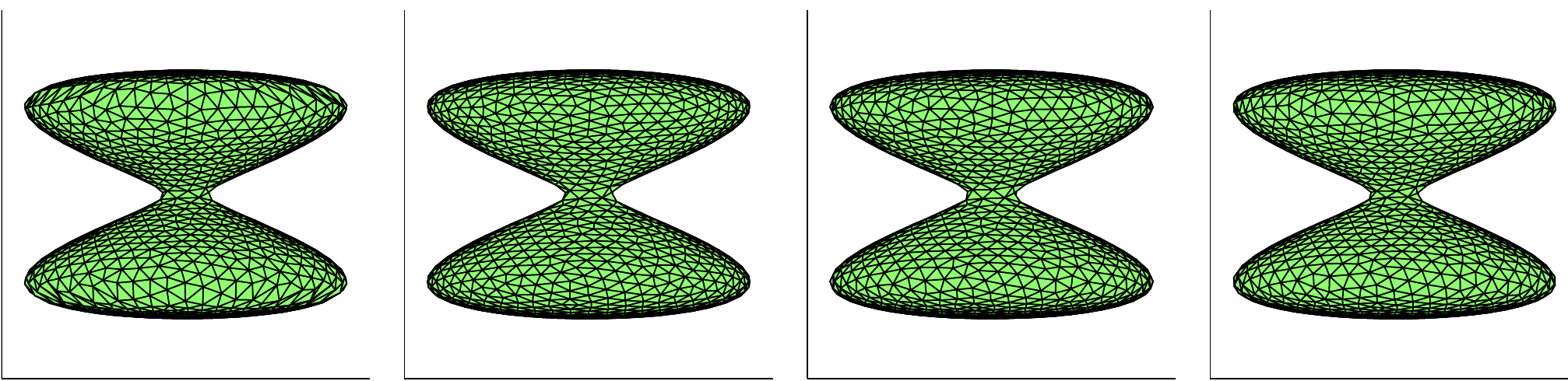}\\
	\includegraphics[width=.99\textwidth]{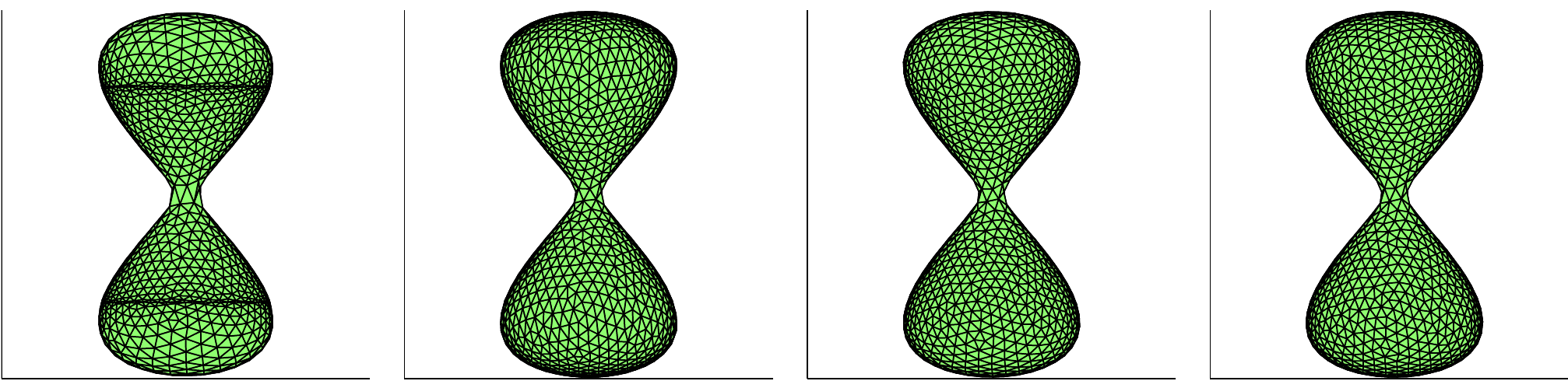}\\
	\caption{Surface evolution of \eqref{eq: surface E-S} using pure normal movement \eqref{eq: normal movement ODE E-S} (first from left); literature ALE map \eqref{eq: ALE map E-S} (second); Runge--Kutta ALE map (third); splitting ALE map (fourth), times $t=0, 0.2, 0.4, 0.6$ (from top to bottom).}
	\label{figure: E-S evolutions}
\end{figure}

\bigskip
To illustrate the good qualities of the DAE model \eqref{eq: constrained spring system} we evolve the surface \eqref{eq: surface E-S} with all four methods.
In Figure~\ref{figure: E-S evolutions} we can observe the evolutions of the discrete initial surface $\Ga_h(0)$ over $[0,0.6]$.
It can be nicely observed that the quality of meshes obtained by both DAE approaches are very similar, however the splitting approach is much faster.

Firstly, plotted on the left-hand side, the purely normal surface evolution obtained by solving the ODE system \eqref{eq: normal movement ODE E-S}. We have used here the explicit Euler method, with step size $\tau=0.001$. Although, this method is clearly not the best choice, we used it in order to illustrate the performance of the proposed ALE algorithms.

Secondly, plotted second from the left, the ALE map \eqref{eq: ALE map E-S} proposed in the literature, cf.\ \cite{ElliottStyles_ALEnumerics}, based on the structure of the distance function \eqref{eq: surface E-S}.

Thirdly, plotted third from the left, the Runge--Kutta solution of the DAE system \eqref{eq: constrained spring system}, using the Radau IIA method with $s=3$ stages, with a stepsize $\tau=0.001$, and $k=500$, $p=0.4$ in \eqref{eq: length function}.

Finally, plotted on the right-hand side, the ALE map obtained by the splitting method. The evolution and the ALE map is computed exactly as described in Section~\ref{section: splitting}. The $v$-system is solved by the explicit Euler method, with $\tau=0.01$, while the $w$-system is solved using the classical Runge--Kutta method of order four with $25$ substeps, and again $p=0.4$

It can be nicely observed visually that both ALE approaches provide meshes of similar quality as the literature ALE map.

\clearpage

%%%%%%%%%%%%%%%%%%%%%%%%%%%%
% mesh quality plots
%%%%%%%%%%%%%%%%%%%%%%%%%%%%
\begin{figure}[ht!]
	\centering
	\includegraphics[width=.99\textwidth]{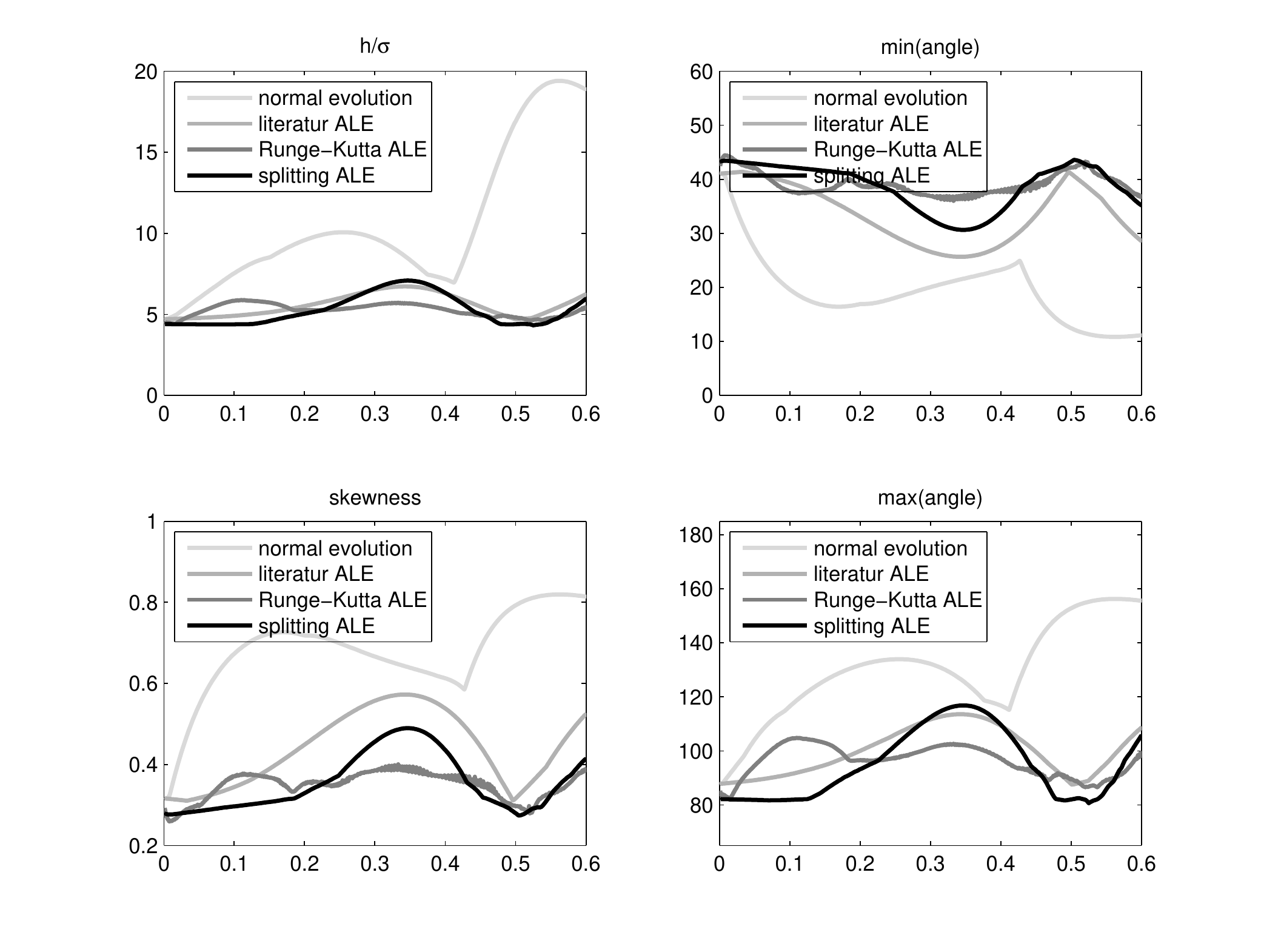}\\
	\caption{Mesh quality measures plotted against time}
	\label{figure: E-S mesh quality}
\end{figure}

%\clearpage

In Figure~\ref{figure: E-S mesh quality} we plotted the evolution of four mesh quality measures (cf.~\cite{Field_meshqualy}) for all four surface evolutions (described above) against time.

In the top left the maximal ratio of element size and the radius of the inscribed circle,
\begin{equation*}
r(t)=\max_{E(t)\in\Th\t} \frac{h_{E(t)}}{\sigma_{E(t)}}
\end{equation*}
can be observed, where $h_E$ is the maximal edge length and $\sigma_E$ is the radius of the inscribed circle of a triangle $E\in\Th\t$. The same mesh quality measure is also used in \cite{ElliottFritz_DT}.
The plots on the right-hand side show minimum and maximum angles ($\alpha_{\textnormal{min}}\t$ and $\alpha_{\textnormal{max}}\t$, top and bottom, respectively) of the mesh, i.e.\
\begin{equation*}
\alpha_{\textnormal{min}}\t = \min_{E\in\Th\t}  \min \alpha_E , \qquad \textnormal{ and } \qquad \alpha_{\textnormal{max}}\t = \max_{E\in\Th\t}  \max \alpha_E ,
\end{equation*}
where $\alpha_E$ contains the three angles of the triangle $E\in\Th\t$.
Finally, the bottom left plot shows the maximal skewness (also called equiangular skew) $s(t)= \max_{E\t\in\Th\t} s_{E\t}$, where $s_E$ is the skewness of a triangle $E\in\Th\t$, and it is defined as
\begin{equation*}
s_E = \max \Big\{ \frac{\max\alpha_E - 60^\circ}{120^\circ}, \frac{60^\circ - \min\alpha_E}{60^\circ} \Big\} \in [0,1].
\end{equation*}
The skewness measures how irregular the triangle $E$ is. For example, for a regular triangle $s_E=0$ and for highly irregular ones $s_E$ tends to one. Usually, a mesh is considered good with skewness less than $0.5$, while non-acceptable if its skewness exceeds $0.8$.

In Figure~\ref{figure: E-S mesh quality} it can be clearly observed that all three ALE maps (literature and DAE solved by Runge--Kutta or splitting method) provide significantly better meshes as the purely normal evolution. Also, both ALE maps from the numerical solutions of the DAE system provide slightly better meshes then the one suggested in the literature.

%\begin{figure}
%    \centering
%    \includegraphics[width=.99\textwidth]{E-S_surmesh_quality_h0.1_0.01}\\
%    \includegraphics[width=.99\textwidth]{E-S_surmesh_quality_h0.1_0.2}\\
%    \caption{Mesh quality measures for \eqref{eq: surface E-S} using pure normal movement \eqref{eq: normal movement ODE E-S}; literature ALE map \eqref{eq: ALE map E-S}; ALE map based on DistMesh; at times $t=0, 0.2$.}
%\end{figure}
%\begin{figure}
%    \centering
%    \includegraphics[width=.99\textwidth]{E-S_surmesh_quality_h0.1_0.4}\\
%    \includegraphics[width=.99\textwidth]{E-S_surmesh_quality_h0.1_0.6}\\
%    \caption{Mesh quality measures for \eqref{eq: surface E-S} using pure normal movement \eqref{eq: normal movement ODE E-S}; literature ALE map \eqref{eq: ALE map E-S}; ALE map based on DistMesh; at times $t=0.4, 0.6$ .}
%\end{figure}

\subsubsection{Surface with four holes \cite{ElliottVenkataraman_ALE}}
%The following example have been used previously in \cite{ElliottVenkataraman_ALE}.
Let the closed surface $\Ga\t$ be given by the zero level set of the distance function
\begin{equation}
\label{eq: surface 4H}
d(x,t) = \frac{x_1^2}{K\t^2} + G(x_2^2) + K\t^2 G\Big( \frac{x_3^2}{L\t^2} \Big)-1, \quad \textrm{i.e.,} \quad \Ga\t = \{ x \in \R^3 \ \big| \ d(x,t) = 0 \}.
\end{equation}
Here the functions $G$, $L$ and $K$ are given by
\begin{align*}
G(s) =&\ 31.25 s (s - 0.36) (s - 0.95) , \\
L(t) =&\ 1 + 0.3\sin(4\pi  t) , \\
K(t) =&\ 0.1 + 0.01\sin(2\pi  t) .
\end{align*}
The normal velocity $v$ describes the surface evolution at the nodes by the ODEs:
\begin{equation}
\label{eq: normal movement ODE 4H}
\diff x_j\t = v(x_j\t,t),
\end{equation}
for $j=1,2,\dotsc,N$, where the surface velocity at $x_j\t$ is again given by the formula \eqref{eq: v def by distance function}.

Finally, the ALE map from the literature is given by
\begin{equation}
\label{eq: ALE map 4H}
(x_i\t)_1= (x_i(0))_1 \frac{K\t}{K(0)}, \quad (x_i\t)_2= (x_i(0))_2, \quad (x_i\t)_3= (x_i(0))_3 \frac{L\t}{L(0)},
\end{equation}
for every $t \in [0,T]$ and for $i=1,2,\dotsc,N$, as suggested in \cite{ElliottVenkataraman_ALE}.

%\begin{figure}
%    \centering
%    \includegraphics[width=.99\textwidth]{4H_ALEmaps_h0.1_0.01}\\
%    \includegraphics[width=.99\textwidth]{4H_ALEmaps_h0.1_0.2}\\
%    \includegraphics[width=.99\textwidth]{4H_ALEmaps_h0.1_0.4}\\
%    \caption{Surface evolution of \eqref{eq: surface 4H} using pure normal movement \eqref{eq: normal movement ODE 4H} (left); literature ALE map \eqref{eq: ALE map 4H} (middle); ALE map based on DistMesh (right); at times $t=0, 0.2, 0.4$ (from top to bottom).}
%\end{figure}
%\begin{figure}
%    \centering
%    \includegraphics[width=.99\textwidth]{4H_ALEmaps_h0.1_0.6}\\
%    \includegraphics[width=.99\textwidth]{4H_ALEmaps_h0.1_0.8}\\
%    \includegraphics[width=.99\textwidth]{4H_ALEmaps_h0.1_1}\\
%    \caption{Surface evolution of \eqref{eq: surface 4H} using pure normal movement \eqref{eq: normal movement ODE 4H}(left); literature ALE map \eqref{eq: ALE map 4H} (middle); ALE map based on DistMesh (right); at times $t=0.6, 0.8, 1$ (from top to bottom).}
%\end{figure}

\begin{figure}[ht!]
	\centering
	\includegraphics[width=.99\textwidth,height=0.166\textheight]{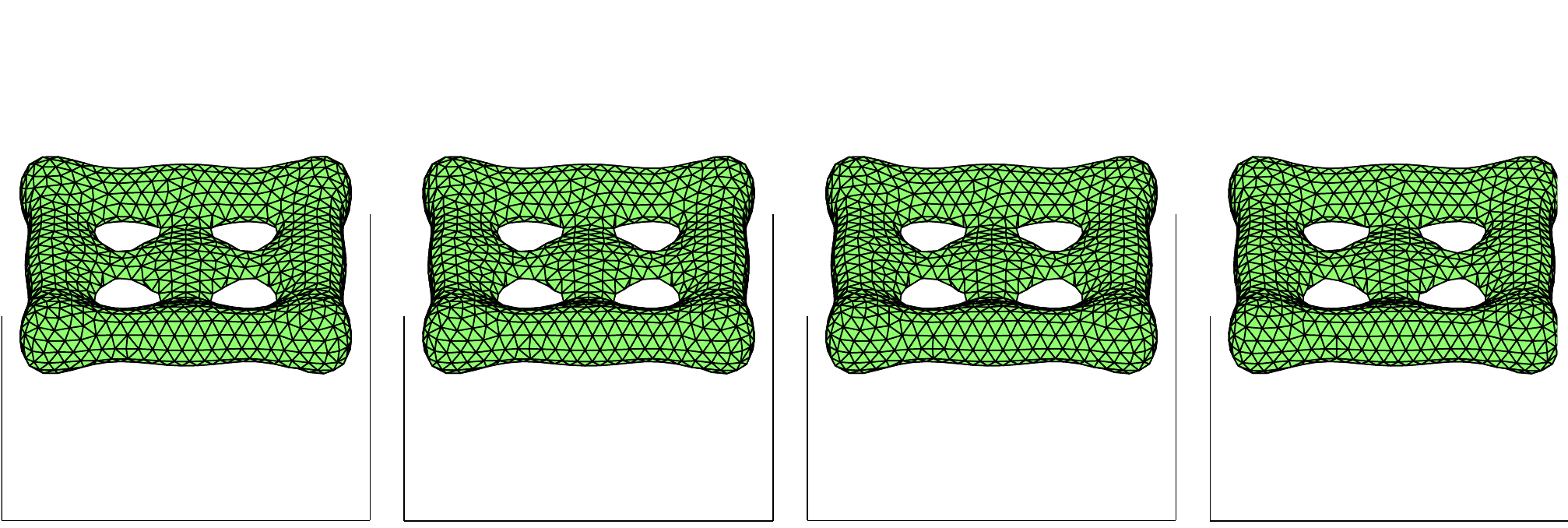}\\
	\includegraphics[width=.99\textwidth,height=0.166\textheight]{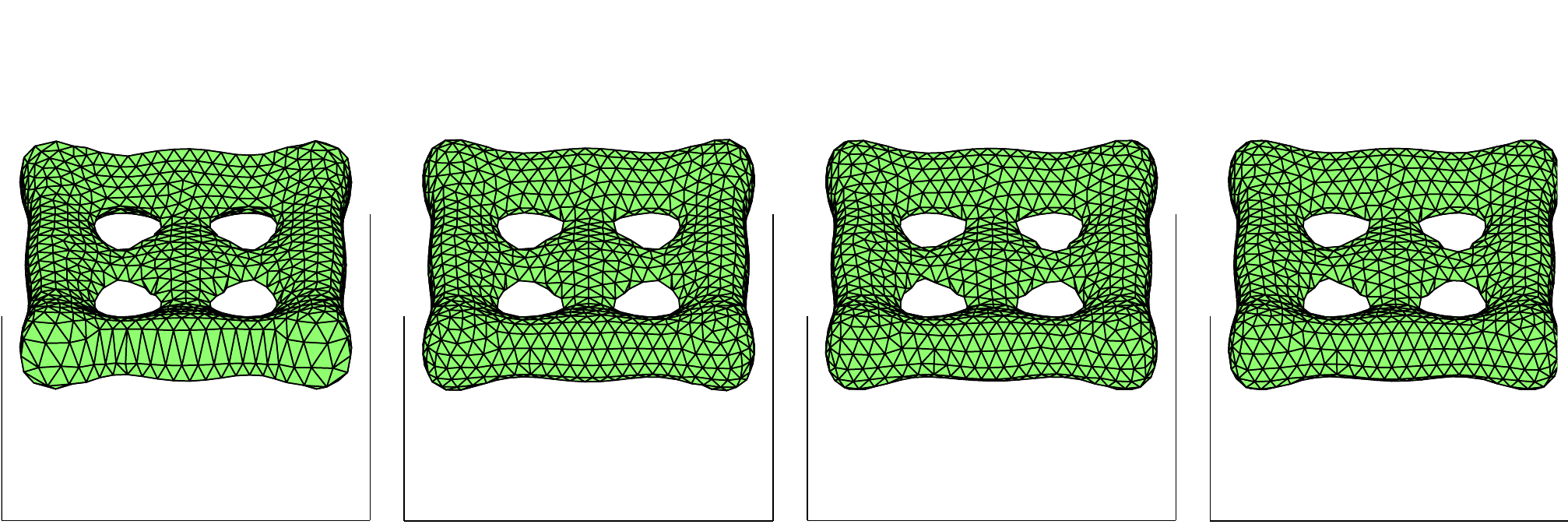}\\
	\includegraphics[width=.99\textwidth,height=0.166\textheight]{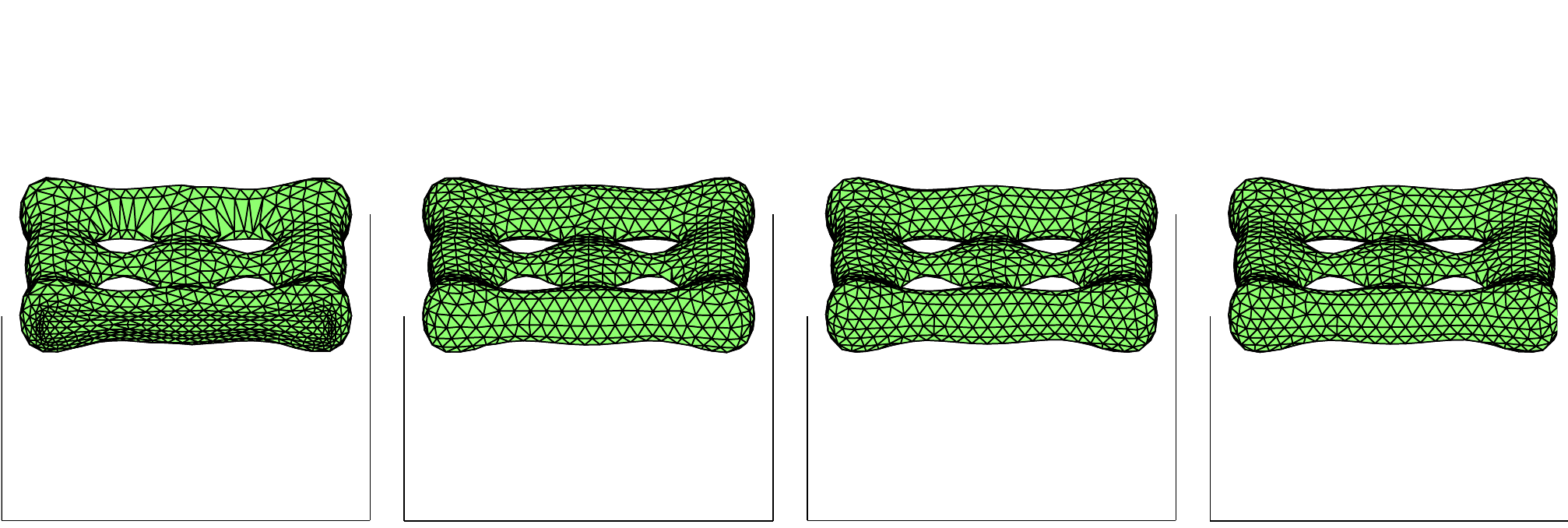}\\
	\includegraphics[width=.99\textwidth,height=0.166\textheight]{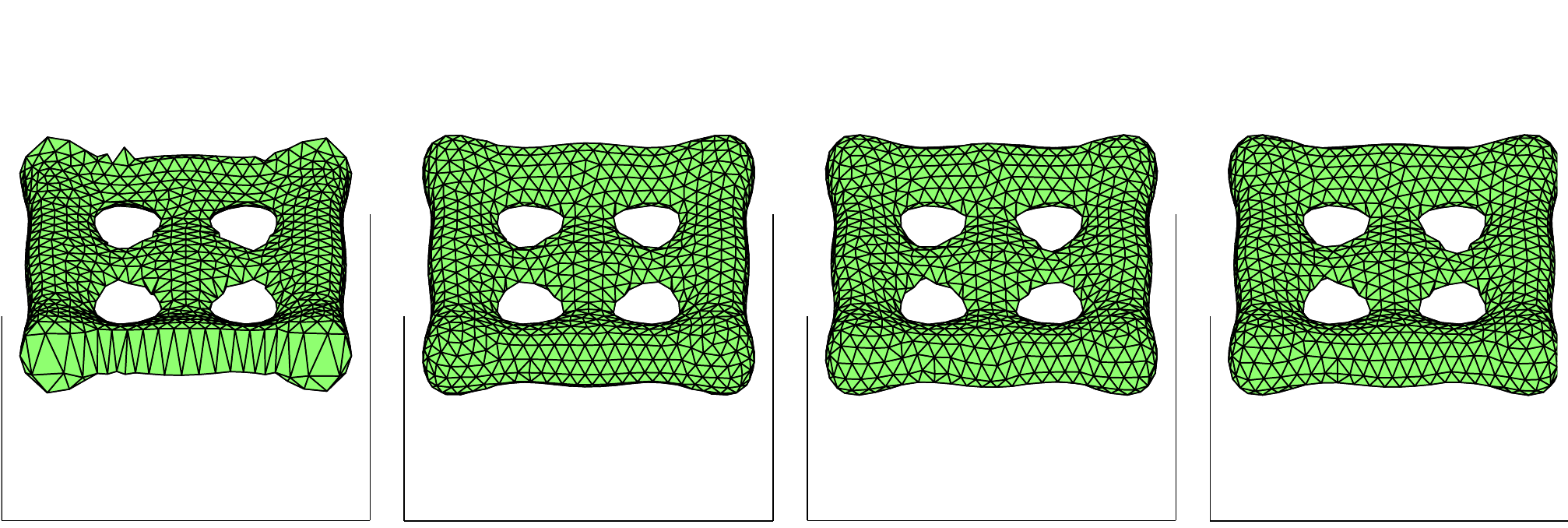}\\
	\includegraphics[width=.99\textwidth,height=0.166\textheight]{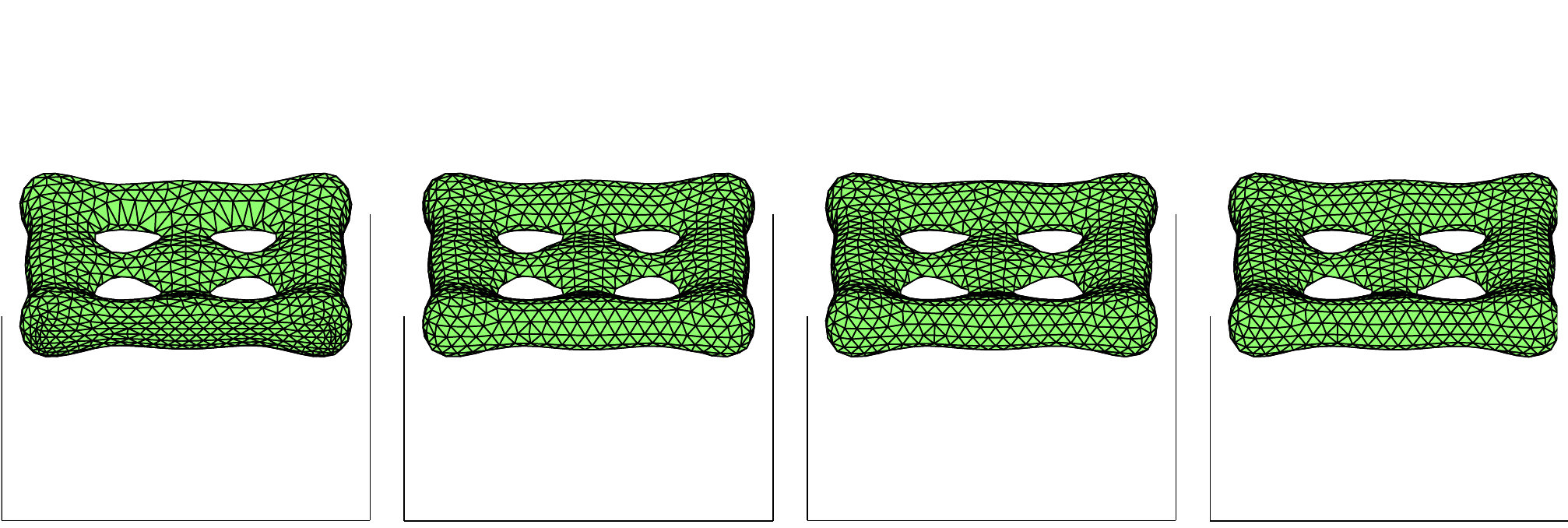}\\
	\includegraphics[width=.99\textwidth,height=0.166\textheight]{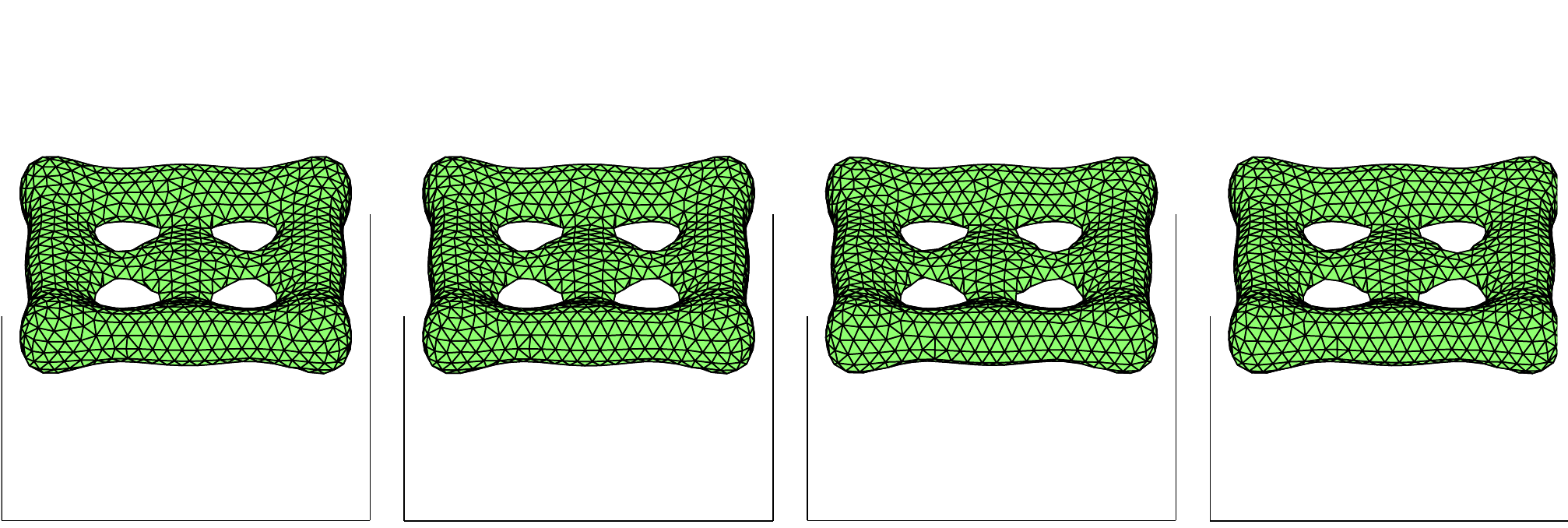}\\
	\caption{Surface evolution of \eqref{eq: surface 4H} using pure normal movement \eqref{eq: normal movement ODE 4H} (first from left); literature ALE map \eqref{eq: ALE map 4H} (second); Runge--Kutta ALE map (third); splitting ALE map (fourth); at times $t=0, 0.2, 0.4, 0.6, 0.8, 1$ (from top to bottom).}
	\label{figure: 4H surface evolution}
\end{figure}

In Figure~\ref{figure: 4H surface evolution}, the surface $\Ga(0)$ (cf.~\eqref{eq: surface 4H}) is evolved over $[0,1]$ again by all four methods, exactly as for Figure~\ref{figure: E-S evolutions}. Again, the poor meshes of the normal movement can be nicely observed, in contrast with the quasi-regular meshes obtained by the ALE maps.

%%%%,height=0.16\textheight
\clearpage

%%%%%%%%%%%%%%%%%%%%%%%%%%%%
% mesh quality plots
%%%%%%%%%%%%%%%%%%%%%%%%%%%%
\begin{figure}[ht!]
	\centering
	\includegraphics[width=.99\textwidth]{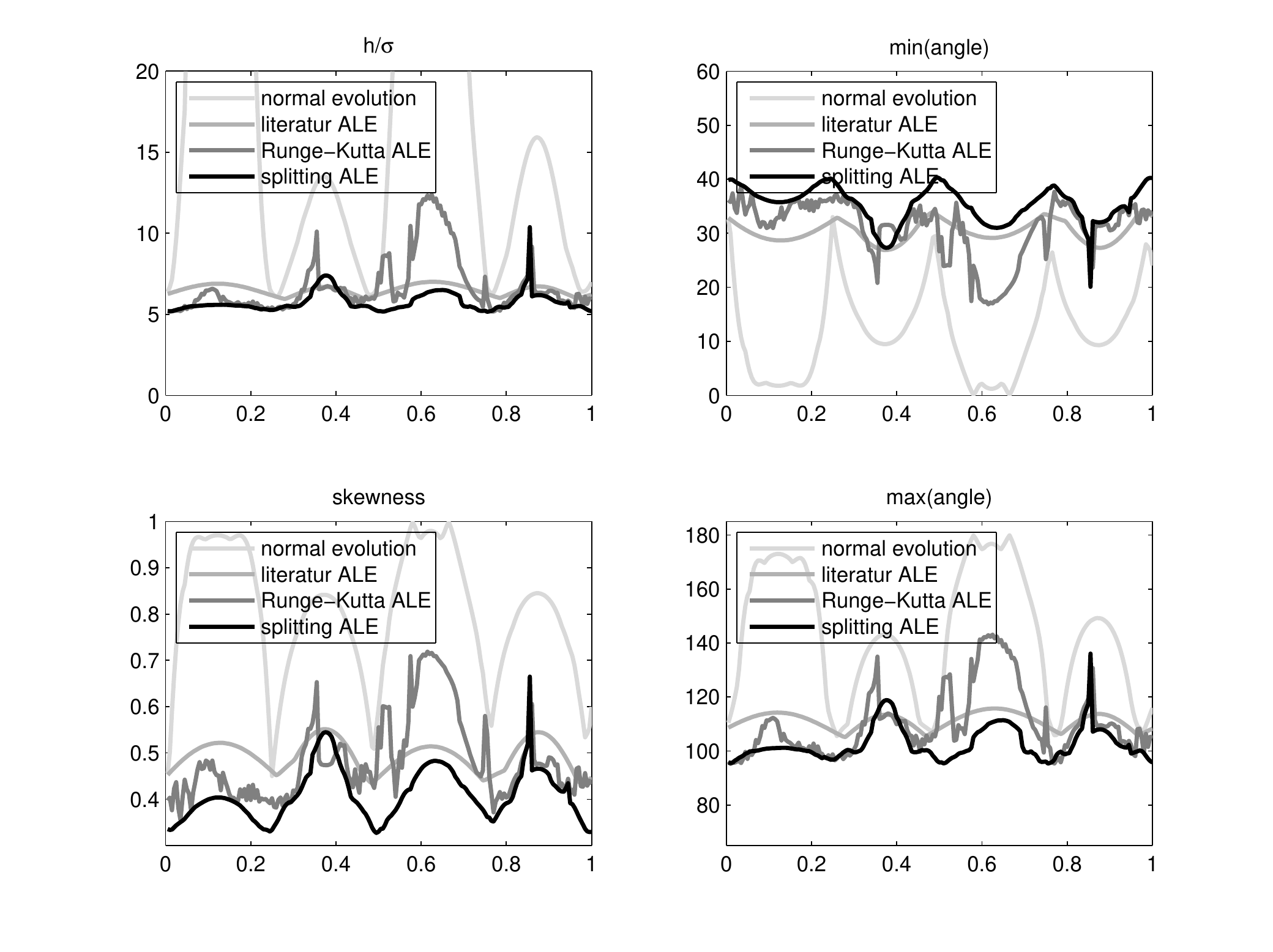}\\
	\caption{Mesh quality measures plotted against time}
	\label{figure: 4H mesh quality}
	
\end{figure}
Similarly, the mesh quality measures are plotted in Figure~\ref{figure: 4H mesh quality}. All three ALE methods provide much better meshes than the normal evolution, while they still yield meshes of similar quality.

%\begin{figure}
%    \centering
%    \includegraphics[width=.99\textwidth]{4H_surmesh_quality_h0.1_0.01}\\
%    \includegraphics[width=.99\textwidth]{4H_surmesh_quality_h0.1_0.2}\\
%    \caption{Mesh quality measures for \eqref{eq: surface 4H} using pure normal movement \eqref{eq: normal movement ODE 4H}; literature ALE map \eqref{eq: ALE map 4H}; ALE map based on DistMesh; at times $t=0, 0.2$ .}
%\end{figure}
%\begin{figure}
%    \centering
%    \includegraphics[width=.99\textwidth]{4H_surmesh_quality_h0.1_0.4}\\
%    \includegraphics[width=.99\textwidth]{4H_surmesh_quality_h0.1_0.6}\\
%    \caption{Mesh quality measures for \eqref{eq: surface 4H} using pure normal movement \eqref{eq: normal movement ODE 4H}; literature ALE map \eqref{eq: ALE map 4H}; ALE map based on DistMesh; at times $t=0.4, 0.6$ .}
%\end{figure}
%\begin{figure}
%    \centering
%    \includegraphics[width=.99\textwidth]{4H_surmesh_quality_h0.1_0.8}\\
%    \includegraphics[width=.99\textwidth]{4H_surmesh_quality_h0.1_1}\\
%    \caption{Mesh quality measures for \eqref{eq: surface 4H} using pure normal movement \eqref{eq: normal movement ODE 4H}; literature ALE map \eqref{eq: ALE map 4H}; ALE map based on DistMesh; at times $t=0.8, 1$ .}
%\end{figure}

\subsection{Error behaviour of evolving surface parabolic PDEs}

We tested the performance of the spring ALE algorithm by using it in the numerical solution of an evolving surface PDE. As a well studied example, we carried out the same experiments over the evolving surface of Section~\ref{section: E-S example}, which have been also used in \cite{ElliottStyles_ALEnumerics,ElliottVenkataraman_ALE,KovacsPower_ALE}.

The evolving surface PDE \eqref{eq: evolving surface PDE} is discretized using ALE ESFEM, described in detail in \cite{ElliottVenkataraman_ALE,KovacsPower_ALE}. The inhomogeneity $f$ is chosen such that the true solution is known: $u(x,t) = e^{-6t} x_1 x_2$.

\subsubsection{Discretization errors}
\begin{figure}[ht!]
	\centering
	\includegraphics[width=.99\textwidth,height=0.3\textheight]{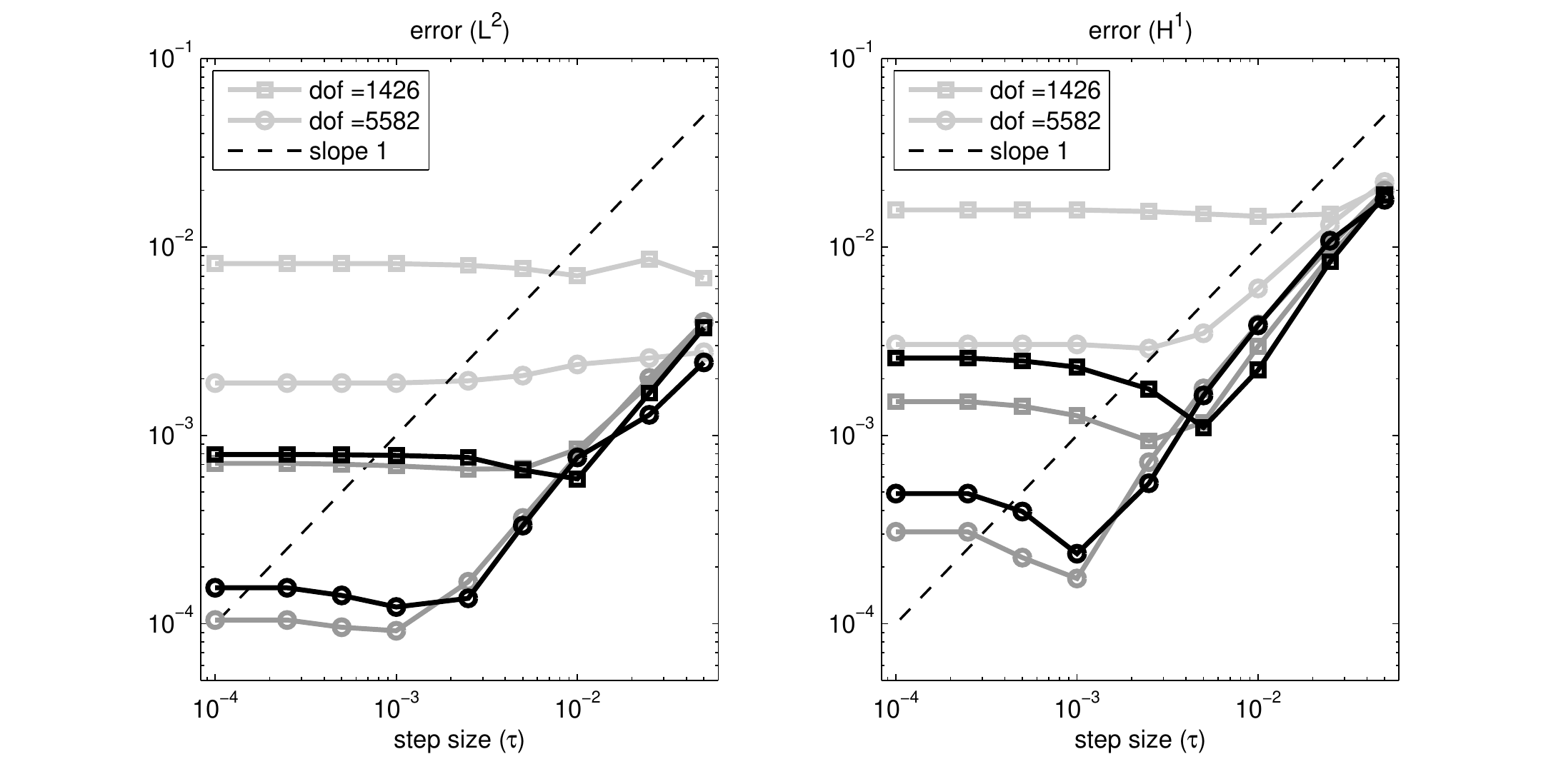}\\
	\caption{$L^2$ and $H^1$ norms of the errors for different spatial refinements plotted against the stepsize for the normal movement, literature and splitting ALE maps (light grey, grey and black, respectively) at time $T=1$}
	\label{figure: conv}
\end{figure}

The errors of the obtained lifted numerical solution is calculated in the following norms at time $T=N\tau=1$:
\begin{equation*}
\|u(\cdot,N \tau)-u_h^N\|_{L^2(\Ga(N\tau))} \andquad \|\nbg(u(\cdot,N \tau)-u_h^N)\|_{L^2(\Ga(N\tau))} .
\end{equation*}
The logarithmic plots, in Figure~\ref{figure: conv}, show the usual convergence plots: the two errors against time step size $\tau$.

%using the following norm and seminorm:
%\begin{alignat*}{3}
%    L^\infty(L^2):&\qquad \max_{1\leq n \leq N}\|u_h^n - u(\cdot,t_n)\|_{L^2(\Ga(t_n))},\\
%    L^2(H^1):&\qquad  \Big(\tau \sum_{n=1}^{N} \|\nb_{\Ga(t_n)}\big(u_h^n - u(\cdot,t_n)\big)\|_{L^2(\Ga(t_n))}^2\Big)^{1/2}.
%\end{alignat*}

Figure~\ref{figure: conv} only serves as an illustration that the meshes from the literature ALE map and from the DAE system yield errors of the same magnitude. More spatial refinements would obfuscate the readability of the plots. For more detailed convergence tests we refer to \cite{ElliottVenkataraman_ALE,KovacsPower_ALE}.

Figure~\ref{figure: conv} shows the errors for all three methods concerning the evolution of the surface (normal evolution, literature ALE and splitting ALE). On the left we show the errors in the $L^2$ norm, while on the right in the $H^1$ norm. Different shades correspond to different ALE maps (light grey, grey and black, respectively), while the lines with different markers are corresponding to different mesh refinements.

As usual we can observe two regions in Figure~\ref{figure: conv}: A region where the time discretization error dominates, matching to the convergence rates of the theoretical results for the backward Euler method (order 1, not the reference line). In the other region, with smaller stepsizes, the spatial discretization error is dominating (the error curves are flatting out).
%The convergence in space and in time can both be nicely observed in agreement with the theoretical results of \cite{ElliottVenkataraman_ALE,KovacsPower_ALE} (note the reference lines).

\subsubsection{Computation times} 
\begin{figure}[htbp]
	\centering
	\includegraphics[width=.99\textwidth]{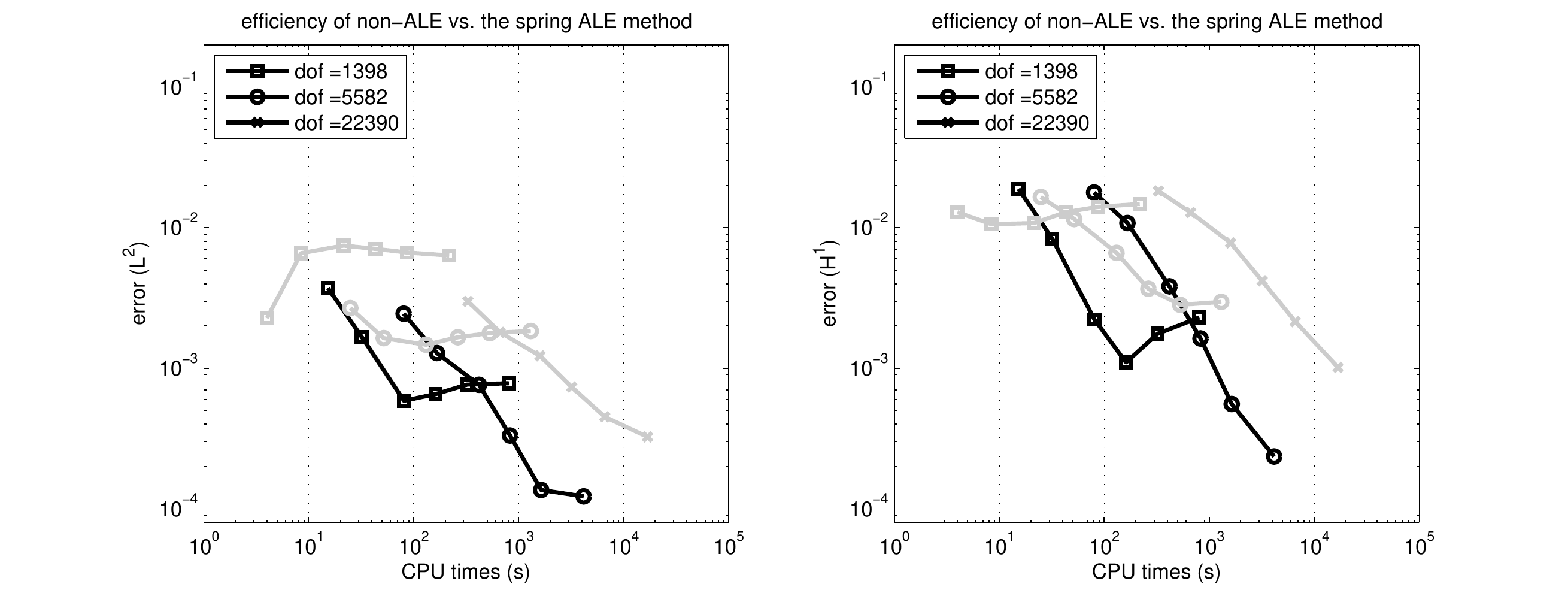}\\
	\caption{$L^2$ and $H^1$ norms of the errors for different spatial refinements and time stepsizes plotted against the CPU times (in seconds) for the normal movement and splitting ALE map (light grey and black, respectively)}
	\label{figure: CPU}
\end{figure}

In Figure~\ref{figure: CPU} we also compare the computational times for the pure normal movement and splitting ALE approach (light grey and black in the figure). In Figure~\ref{figure: CPU} the errors at time $T=1$ (both in the $L^2$ and $H^1$ norms, left and right, respectively) are plotted against the CPU times for different spatial refinements (denoted by different markers: {\footnotesize $\square$}, $\circ$, $\times$) and different stepsizes (which are $\tau=0.05, 0.025, 0.01, 0.005$, $0.0025$, $0.001$, corresponding to the markers). The CPU times include all computations: normal evolution of the surface, surface matrix assemblies, solution of the linear system, and also the computation of the ALE mesh in the ALE case. A computational trade-off can be clearly observed, (see, for instance the two rightmost lines in each plot), that the same errors can be achieved by the ALE method on a much coarser mesh (having only quarter as many nodes as the non-ALE mesh), at a computational cost reduction of factor $10$. Again, in the region with smaller stepsizes, the spatial discretization error is dominating (hence the curves are flatting out).

\subsection{Meshes with angle conditions}

We also report on a somewhat simplified version of the force function proposed in Section~\ref{section: acute meshes} in order to construct acute or nonobtuse meshes. 

\begin{figure}[htbp]
	\centering
	\includegraphics[width=1\textwidth]{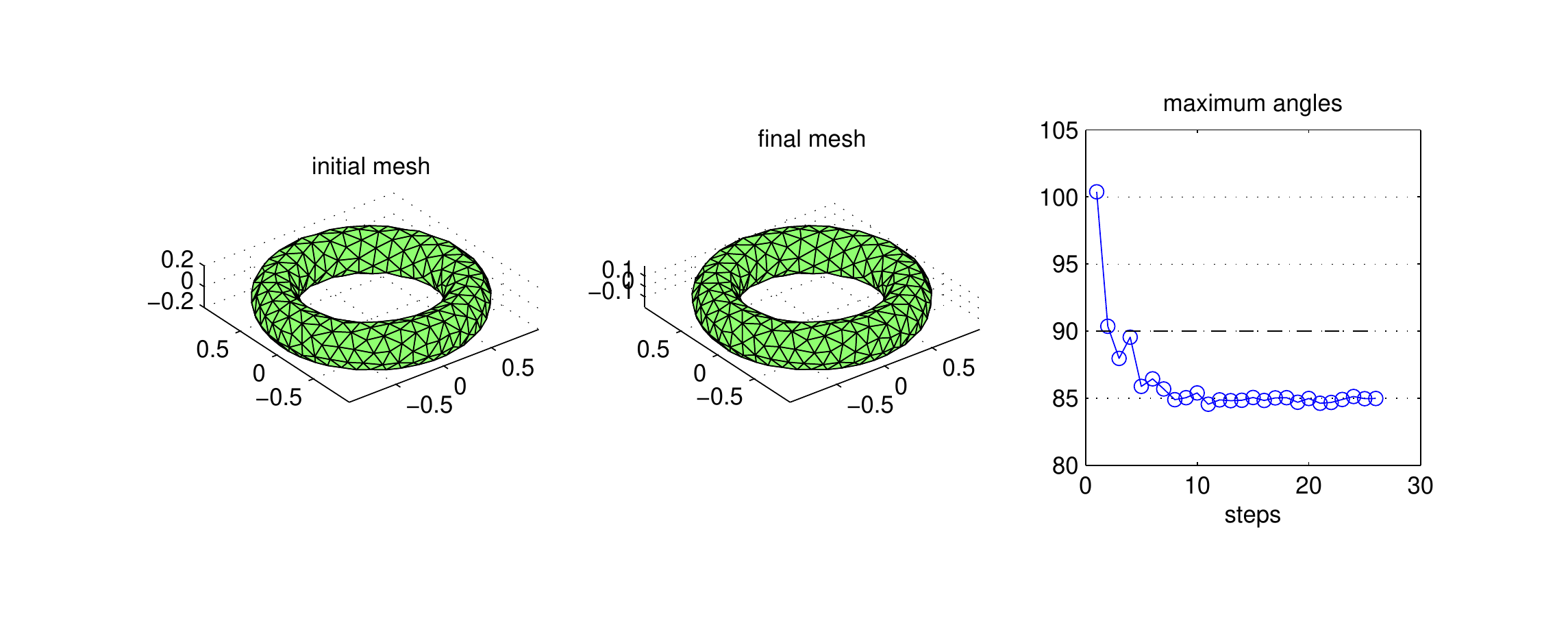}\\
	\caption{Acute triangulation of a torus}
	\label{figure: torus}
\end{figure}

We consider here a stationary example of a torus. We constructed a triangulation using DistMesh \cite{distmesh}, seen left in Figure~\ref{figure: torus}, which is not acute, the largest angle of the mesh is above $100^\circ$ (see the first entry in rightmost graph in Figure~\ref{figure: torus}). 

As the surface is stationary here, we have the modified DAE system
\begin{equation*}
\begin{aligned}
\diff \bfx\t =&\  k F(\bfx\t) + k_\alpha F_\alpha(\bfx\t) - D(\bfx\t)^T \bflambda\t , \\
d(\bfx\t,t) =&\ 0 ,
\end{aligned}
\end{equation*}
with an additional force function $F_\alpha$, which is used to eliminate non-acute angles. It is defined analogously as the function $F$, but using a length function $f_\alpha(e)$ based on the angle $\alpha_e$ opposite to the edge $e$, instead of the formula in \eqref{eq: length function}. If the angle $\alpha_e$ is larger than a user given tolerance, the desired length of $e$ is set to a value based on the law of cosines. As the important force is now $F_\alpha$ we set its spring constant $k_\alpha=4k$. The original force function $F$ (with $p=0.1$ in \eqref{eq: length function}) is only kept for smoothing reasons. We note here, that these parameters are not claimed to be used universally.

The numerical solution of the above DAE system yields meshes with favourable angle properties. We set here the angle tolerance to $85^\circ$ and integrate the system over 25 steps. The maximum angle in each time step can be seen on the right in Figure~\ref{figure: torus}, it can be observed that the angles quickly drop around the desired value, while the finally obtained acute mesh can be seen in the middle of Figure~\ref{figure: torus}.

%\clearpage
\section*{Acknowledgement}
We thank Christian Lubich for our many discussions on differential algebraic problems, and for his invaluable comments and remarks. We thank Hanna Walach for her helpful suggestions regarding the manuscript. We also thank Christian Power for our early discussions on the topic.

This work was supported by the Deutsche Forschungsgemeinschaft (DFG) through SFB 1173.

%\pagebreak

\bibliographystyle{siamplain}
%\bibliography{evolving_surface_literature}

\begin{thebibliography}{10}
	
	\bibitem{DUNE}
	{\sc M.~Blatt, A.~Burchardt, A.~Dedner, C.~Engwer, J.~Fahlke, B.~Flemisch,
		C.~Gersbacher, C.~Gr{\"a}ser, F.~Gruber, C.~Gr{\"u}ninger, et~al.}, {\em The
		distributed and unified numerics environment, version 2.4}, Archive of
	Numerical Software, 4 (2016), pp.~13--29.
	
	\bibitem{Bonito2013a}
	{\sc A.~Bonito, I.~Kyza, and R.~H. Nochetto}, {\em Time--discrete higher--order
		{ALE} formulations: a priori error analysis}, Numer. Math.,  (2013),
	pp.~577--604.%, \url{https://doi.org/10.1007/s00211-013-0549-3}.
	
	\bibitem{Bonito2013}
	{\sc A.~Bonito, I.~Kyza, and R.~H. Nochetto}, {\em Time--discrete higher--order
		{ALE} formulations: stability}, SIAM J. Numer. Anal., 51 (2013),
	pp.~577--604.%, \url{https://doi.org/10.1137/120862715}.
	
	\bibitem{Chang}
	{\sc W.~Chang}, {\em Surface reconstruction from points}, UCSD CSE Technical
	Report, CS2008-0922,  (2007), pp.~1--18.
	
	\bibitem{Dziuk88}
	{\sc G.~Dziuk}, {\em Finite elements for the {B}eltrami operator on arbitrary
		surfaces}, Partial differential equations and calculus of variationsLecture Notes in Math., 1357, Springer, Berlin, (1988), pp.~142--155.
	
	\bibitem{DziukMCF}
	{\sc G.~Dziuk}, {\em An algorithm for evolutionary surfaces}, Numerische Mathematik, 58(1) (1990), pp.~603--611.
	
	\bibitem{DziukElliott_ESFEM}
	{\sc G.~Dziuk and C.~M.~Elliott}, {\em Finite elements on evolving surfaces}, IMA
	Journal of Numerical Analysis, 27 (2007), pp.~262--292.%, \url{https://doi.org/10.1093/imanum/drl023}.
	
	\bibitem{DziukElliott_acta}
	{\sc G.~Dziuk and C.~M.~Elliott}, {\em Finite element methods for surface {PDE}s},
	Acta Numerica, 22 (2013), pp.~289--396.%,  \url{https://doi.org/10.1017/S0962492913000056}.
	
	\bibitem{DziukElliott_L2}
	{\sc G.~Dziuk and C.~M.~Elliott}, {\em {$L^2$}--estimates for the evolving surface
		finite element method}, Mathematics of Computation, 82 (2013), pp.~1--24.
	
	\bibitem{DziukKronerMuller}
	G.~Dziuk, D.~Kr\"{o}ner, and T.~M\"{u}ller, {\em Scalar conservation laws on moving hypersurfaces}, Interfaces and Free Boundaries, 15(2):203--236, 2013.
	
	\bibitem{ElliottFritz_DT}
	{\sc C.~M.~Elliott and H.~Fritz}, {\em On algorithms with good mesh properties for
		problems with moving boundaries based on the harmonic map heat flow and the
		{DeTurck} trick}, SMAI-Journal of Computational Mathematics, 2:141--176, (2016).
	
	\bibitem{ElliottFritz}
	{\sc C.~M.~Elliott and H.~Fritz}, {\em On approximations of the curve shortening flow and of the mean curvature flow based on the {DeTurck} trick}, IMA Journal of Numerical Analysis, doi: 10.1093/imanum/drw020, (arXiv:1602.07143),  (2016).
	
	\bibitem{ElliottStyles_ALEnumerics}
	{\sc C.~M.~Elliott and V.~Styles}, {\em An {ALE} {ESFEM} for solving {PDE}s on evolving surfaces}, Milan Journal of Mathematics, 80 (2012), pp.~469--501.%,
	%  \url{https://doi.org/10.1007/s00032-012-0195-6}.
	
	\bibitem{ElliottVenkataraman_ALE}
	{\sc C.~M. Elliott and C.~Venkataraman}, {\em {E}rror analysis for an {ALE}
		evolving surface finite element method}, Numerical Methods for Partial
	Differential Equations, 31 (2015), pp.~459--499.
	
	\bibitem{FarhatLesoinneMaman}
	{\sc C.~Farhat, M.~Lesoinne, and N.~Maman}, {\em Mixed explicit-implicit time
		integration of coupled aeroelastic problems: {T}hree-field formulation,
		geometric conservation and distributed solution}, International Journal for
	Numerical Methods in Fluids, 21 (1995), pp.~807--835.
	
	\bibitem{Field_meshqualy}
	{\sc D.~A.~Field}, {\em Qualitative measures for initial meshes}, International Journal for Numerical Methods in Engineering, 47(4) (2000), pp.~887--906.
	
	\bibitem{FormaggiaNobile1999}
	{\sc L.~Formaggia and F.~Nobile}, {\em A stability analysis for the arbitrary
		{L}agrangian {E}ulerian formulation with finite elements}, East-West J.
	Numer. Math., 7 (1999), pp.~105--131.
	
	\bibitem{FormaggiaNobile2004}
	{\sc L.~Formaggia and F.~Nobile}, {\em Stability analysis of second-order time
		accurate schemes for {ALE--FEM}}, Computer methods in applied mechanics and
	engineering, 193 (2004), pp.~4097--4116.
	
	\bibitem{lumped}
	{\sc M.~Frittelli, A.~Madzvamuse, I.~Sgura, and C.~Venkataraman}, {\em Lumped
		finite element method for reaction-diffusion systems on compact surfaces},
	arXiv:1609.02741,  (2016).
	
	\bibitem{HairerLubichRoche}
	{\sc E.~Hairer, {Ch}.~Lubich, and M.~Roche}, {\em The numerical solution of differential-algebraic systems by Runge-Kutta methods}. Vol.~1409. Springer, Berlin (2006).
	
	\bibitem{HairerWannerII}
	{\sc E.~Hairer and G.~Wanner}, {\em Solving Ordinary Differential Equations
		II.: Stiff and differetial--algebraic problems}, Springer, Berlin {S}econd~ed.,
	(1996).
	
	\bibitem{Hoppe}
	{\sc H.~Hoppe, T.~DeRose, T.~Duchamp, J.~McDonald, and W.~Stuetzle}, {\em
		Surface reconstruction from unorganized points}, Proceedings of the 19th
	annual conference on Computer Graphics and interactive techniques, (1992).
	
	\bibitem{surfDMP}
	{\sc J.~Kar{\'a}tson, B.~Kov{\'a}cs, and S.~Korotov}, {\em Discrete maximum
		principles for nonlinear elliptic finite element problems on {R}iemannian
		manifolds with boundary}, arXiv:1701.00424, (2016).
	
	\bibitem{soldriven}
	{\sc B.~Kov\'{a}cs, B.~Li, {Ch}.~Lubich, and C.~{Power Guerra}}, {\em Convergence of finite elements on an evolving surface driven by diffusion on the surface}, to appear in {\em Numerische Mathematik}, arXiv:1607.07170,
	(2016).
	
	\bibitem{KovacsPower_ALE}
	{\sc B.~Kov\'{a}cs and C.~{Power Guerra}}, {\em Higher--order time
		discretizations with {ALE} finite elements for parabolic problems on evolving
		surfaces}, to appear in {\em IMA Journal of Numerical Analysis}, doi: 10.1093/imanum/drw074, arXiv:410.0486,  (2014).
	
	\bibitem{LubichMansourVenkataraman_bdsurf}
	{Ch}.~Lubich, D.E.~Mansour, and C.~Venkataraman.
	\newblock Backward difference time discretization of parabolic differential
	equations on evolving surfaces.
	\newblock {\em IMA Journal of Numerical Analysis}, 33(4):1365--1385, (2013).
	
	\bibitem{diss_Mansour}
	{\sc D.~Mansour}, {\em Numerical Analysis of Partial Differential Equations on
		Evolving Surfaces}, PhD thesis, Universit\"{a}t T\"{u}bingen, (2013).
	\newblock {http://hdl.handle.net/10900/49925}.
	
	\bibitem{distmesh}
	{\sc P.-O. Persson and G.~Strang}, {\em A simple mesh generator in {MATLAB}},
	SIAM Review, 46 (2004), pp.~329--345.
	
\end{thebibliography}

\end{document}